\documentclass[aop,MSNbibl,seceqn,citesort,dvips]{arximspdf}

\doi{10.1214/11-AOP733} %kopijuoti is PTS
\volume{41}
\issue{1}
\pubyear{2013}
\firstpage{170}
\lastpage{205}

\makeatletter
\newtheorem{tm}{Theorem}[section]
\newtheorem{lm}[tm]{Lemma}
\newtheorem{co}[tm]{Corollary}
\newproclaim{re}[tm]{Remark}
\newtheorem{pr}[tm]{Proposition}
\newcommand{\iint}{\int\!\!\!\int}
\makeatother

\begin{document}
\begin{frontmatter}

\title{Laplace approximation for rough differential equation driven by fractional Brownian motion}
\runtitle{Laplace approximation for RDE driven by fBM}

\begin{aug}
\author[A]{\fnms{Yuzuru} \snm{Inahama}\corref{}\ead[label=e1]{inahama@math.nagoya-u.ac.jp}}
\runauthor{Y. Inahama}
\affiliation{Nagoya University}
\address[A]{Graduate School of Mathematics\\
Nagoya University\\
Furocho, Chikusa-ku, Nagoya 464-8602\\
Japan\\
\printead{e1}} %adresu isvedimo komanda gale!
\end{aug}

% HISTORY:
\received{\smonth{9} \syear{2009}}
\revised{\smonth{6} \syear{2011}}

% ABSTRACT
%
\begin{abstract}
We consider a rough differential equation indexed by a small parameter
$\varepsilon>0$.
When the rough differential equation is driven by
fractional Brownian motion with Hurst parameter $H$ ($1/4 <H<1/2$),
we prove the Laplace-type asymptotics for the solution as the parameter
$\varepsilon$ tends to zero.
\end{abstract}

% KEYWORDS
%

\end{frontmatter}

%s1 #&#
\section{Introduction}\label{sec1}
The rough path theory was invented by T. Lyons in~\cite{ly} and
summarized in a book~\cite{lq} with Z. Qian.
See also~\cite{lej,fvbook,lcl}.
Roughly speaking, a~rough path is a path coupled with its iterated integrals.
T. Lyons generalized the line integral of one-form along a path to the
one along a rough path.
This is a pathwise integral theory and no probability measure is involved.
In a natural way, an ordinary differential equation (ODE) is generalized.
This is called a rough differential equation (RDE) in this paper.
The corresponding It\^o map is not only everywhere defined, but is also
locally Lipschitz continuous
with respect to the topology of geometric rough path space (Lyons'
continuity theorem).
If a Wiener-like measure is given on the geometric rough path space
or, in other words, if a Brownian rough path is mapped by the It\^o map,
then the solution of the corresponding stochastic differential equation
(SDE) of Stratonovich-type is recovered via rough paths.
In order to investigate the Brownian motion, one only needs the double
integral (i.e., the second level path)
as well as the path itself (i.e., the first level path).
In short, we can obtain the solution of an SDE as the image of a
continuous map.
This is basically impossible in the framework of the usual stochastic calculus.
Recall that, in the usual stochastic calculus,
stochastic integrals and SDEs are defined by the martingale integration theory,
which is quite probabilistic by definition. Therefore, those objects
have no pathwise meaning.

Brownian motion and Brownian rough path are most important and were
studied extensively.
There may be other stochastic processes (i.e., probability measures on
the usual path space), however,
which can be lifted
to probability measures on the geometric rough path space.\vspace*{1pt}
The most typical example is the $d$-dimensional fractional Brownian
motion (fBM) $(w^H_t )_{0 \le t \le1} = (w_t^{H,1} , \ldots,
w_t^{H,d})_{0 \le t \le1} $
with
Hurst parameter $H \in(1/4,1/2]$ (see Coutin--Qian~\cite{cq}).
Recall that, when $H=1/2$, it is the Brownian motion.
It is worth noting that, if $H \in(1/4, 1/3]$, the third level path
plays a role,
unlike the Brownian motion case.
The Schilder-type large deviation for the lift of scaled fBM
was proved by Millet and Sanz-Sole~\cite{mss}.
Combined with Lyons' continuity theorem and the contraction principle,
this fact implies that
the solution of an RDE driven by the lift of scaled fBM
also satisfies large deviation.

According to~\cite{cou,bhoz}, there are several types of path integrals
along fBM, namely, (1) deterministic or pathwise integral, (2) integral
with generalized covariation,
(3)~the divergence operator in the sense of the Malliavin calculus, and
(4) White noise approach.
Clearly, the rough path approach belongs to the first category.

More precisely, we consider the following RDE: for $\varepsilon>0$,
%
%e1.1 #&#
\begin{equation}\label{defrdeintro}
dY^{\varepsilon}_t = \sigma(Y^{\varepsilon}_t ) \varepsilon \,dW^H_t
+\beta(\varepsilon, Y^{\varepsilon}_t) \,dt,\qquad
  Y^{\varepsilon}_0=0.
\end{equation}
Here, $W^H$ is the fractional Brownian rough path (fBRP),
that is, the lift of fBM $w^H$
and $\sigma\in C_b^{\infty} ( {\mathbf{R}}^n , \operatorname{Mat}(n,d))$ and $\beta
\in C_b^{\infty} ([0,1] \times{\mathbf{R}}^n, {\mathbf{R}}^n)$.
Note that $C_b^{\infty}$ denotes the set of bounded smooth functions
with bounded derivatives.

The main purpose of this paper is to prove the Laplace approximation
for (the first level path of) $Y^{\varepsilon}$ as $\varepsilon
\searrow0$.
The precise statement is in Theorem~\ref{thmmain} below.
Apparently, in none of the integrals (1)--(4) has the Laplace
approximation been proved
for the solution of SDE (or RDE)
driven by the scaled fBM.
Note that it is the precise asymptotics of the large deviation.
In this paper,
we will prove it in the framework of the rough path theory for $H \in
(1/4, 1/2)$.
[The case $H >1/2$ is not so interesting from a viewpoint of rough path
analysis.
Our method does not work for the case $H \le1/4$ since the relation
for the Young integral $1/p + 1/q >1$ in equation (\ref{ineqparapq}) below
fails to hold.]

The history of this kind of problem is long. A partial list could be as follows.
First, Azencott~\cite{az} showed this kind of asymptotics for finite
dimensional SDEs,
which is followed by Ben Arous~\cite{ba}.
There are similar results for infinite dimensional SDEs (e.g.,
Albeverio--R\"ockle--Steblovskaya~\cite{ars})
as well as SPDEs (e.g., Rovira--Tindel~\cite{rt}).
In the framework of the Malliavin calculus, there are deep results
on the asymptotics of
the generalized expectation of generalized Wiener functionals
(Takanobu--Watanabe~\cite{tw}, Kusuoka--Stroock~\cite{ks1,ks2},
Kusuoka--Osajima~\cite{ko})
which have applications to the asymptotics for the heat kernels on
Riemannian manifolds.

In the framework of the rough path theory,
Aida studied this problem for finite dimensional Brownian rough paths
and gave a new proof for the results in~\cite{az,ba}.
The same problem for infinite dimensional Brownian rough paths
was studied in~\cite{ina,ik}, which has an application to Brownian
motion over loop groups.

The organization of this paper is as follows:
In Section~\ref{sec2} we give a precise statement of our main result.
In Section~\ref{sec3} we review the rough path theory and fractional Brownian
rough path.
In Section~\ref{sec4} we prove the Hilbert--Schmidt property of the Hessian of the
It\^o map restricted on the Cameron--Martin space $\mathcal{H}^H$ of fBM.
For those who understand the proof of Laplace approximation for
Brownian rough path as in~\cite{aida,ina,ik},
this is the most difficult part,
because the Cameron--Martin space of fBM is not understood very well.
However, thanks to Friz--Victoir's result (Proposition~\ref{prFV}),
such Cameron--Martin paths are Young integrable and, therefore, the
Hessian is computable.
In Section~\ref{sec5} we give a probabilistic representation of (the stochastic
extension of) the Hessian.
In Section~\ref{secproof} we give a proof of the main theorem.
In Section~\ref{sec7} we consider the Laplace approximation for an RDE, which
involves a fractional order term of $\varepsilon>0$.
This has an application to the short time asymptotics of integral
quantities of the solution of a fixed RDE driven by fBM.
(Similar problems were studied in~\cite{bc,nnrt}).

%re1.1 #&#
\begin{re}
All the results in this paper hold for the case $H=1/2$, too, with
trivial modifications.
The only reason we do not treat
the case $H=1/2$ (i.e., the usual Brownian case) is because those
results are already well known in that case.
\end{re}

%%%%%%%%%%%%%%%%%%%%%%%%%%%%%%%%%%%%%%%%%%%%%%%%%%%%%%%%%%%%%%%%%%%%%%
%
%%%%%%%%%%%%%%%%%%%%%%%%%%%%%%%%%%%%%%%%%%%%%%%%%%%%%%%%%%%%%%%%%%%%%%
%
%%%%%%%%%%%%%%%%%%%%%%%%%%%%%%%%%%%%%%%%%%%%%%%%%%%%%%%%%%%%%%%%%%%%%%
%
%s2 #&#
\section{Statement of main result}\label{sec2}
%s2.1 #&#
\subsection{Assumption and main result}\label{sec21}
In this section we state our main results in this paper.
Throughout this paper, the time interval is $[0,1]$ except otherwise stated.
Let $1/4 <H<1/2$
and let $\mathcal{H}^H$ be the Cameron--Martin subspace of the
$d$-dimensional fBM $(w^H_t)_{0 \le t \le1}$.
By Friz--Victoir's result,\vspace*{1pt} which will be explained in Proposition \ref
{prFV} below,
$k \in\mathcal{H}^H$ is of finite $q$-variation for any $ (H+1/2)^{-1} <q<2$.
Hence, the following ODE makes sense in the $q$-variational setting in
the sense of the Young integration:
\[
dy_t = \sigma(y_t) \,dk_t +\beta(0,y_t) \,dt, \qquad  y_0=0.
\]
Note that $y$ is again of finite $q$-variation and we will write
$y=\Psi(k)$.

Now we set the following assumptions.
In short, we assume that there is only one point that attains the
minimum of $F_{\Lambda}$ and
the Hessian at the point is nondegenerate.
These are typical assumptions for Laplace's method of this kind.
The space of continuous paths in ${\mathbf{R}}^n$ with finite $p'$-variation
starting at $0$\vadjust{\goodbreak}
is denoted by $C^{p' \mbox{-}\mathrm{var}}_0 ({\mathbf{R}}^n)$.
Note that the self-adjoint operator $A$ in the fourth assumption turns
out to be Hilbert--Schmidt in Theorem~\ref{thmHSmain} below.

\begin{longlist}[(H1):]
\item[(H1):]
$F$ and $G$ are real-valued bounded continuous functions on $C^{p'
\mbox{-}\mathrm{var}}_0 ({\mathbf{R}}^n)$
for some $p' >1/H$.

\item[(H2):]
The function $F_{\Lambda} := F \circ\Psi+ \|  \cdot  \|^2_{\mathcal{H}^H } /2$
attains its minimum at a unique point $\gamma\in\mathcal{H}^H$.
We will write $\phi^0 = \Psi(\gamma)$.

\item[(H3):]
$F$ and $G$ are $m+3$ and $m+1$ times Fr\'echet differentiable on a
neighborhood $U(\phi^0)$ of $\phi^0 \in C^{p' \mbox{-}\mathrm{var}}_0 ({\mathbf{R}}^n)$,
respectively.
Moreover, there are positive constants $M_1, M_2, \ldots$ such that
\begin{eqnarray*}
| \nabla^j F (\eta) \langle z,\ldots, z \rangle| &\le& M_j \| z \|^j_{p'
\mbox{-}\mathrm{var}} \qquad  (j=1,\ldots, m+3),
\\
| \nabla^j G (\eta) \langle z,\ldots, z \rangle| &\le& M_j \| z \|^j_{p'
\mbox{-}\mathrm{var}} \qquad  (j=1,\ldots, m+1)
\end{eqnarray*}
hold for any $\eta\in U(\phi^0)$ and $z \in C^{p' \mbox{-}\mathrm{var}}_0 ({\mathbf{R}}^n)$.

\item[(H4):]
At the point $\gamma\in\mathcal{H}^H$,
the bounded self-adjoint operator $A$ on $\mathcal{H}^H$, which
corresponds to
the Hessian
$
\nabla^2 (F \circ\Psi) (\gamma) \vert_{\mathcal{H}^H \times\mathcal{H}^H}
$,
is strictly larger than $- \mathrm{Id}_{\mathcal{H}^H}$ (in the form sense).
\end{longlist}

Under these assumptions, the following Laplace-type asymptotics hold.
Explicitly, the constant $c= \nabla F (\phi^0) \langle\theta
^1\rangle$, where
$\theta^1$ will be given in (\ref{deftheta1}) below.
[Below, $Y^{\varepsilon, 1}=(Y^{\varepsilon})^1$ denotes the first
level path of
$Y^{\varepsilon}$.]
%
%th2.1 #&#
\begin{tm}
\label{thmmain}
Let the coefficients $\sigma\dvtx  {\mathbf{R}}^n \to\operatorname{Mat} (n,d)$ and
$\beta
\dvtx  [0,1] \times{\mathbf{R}}^n \to{\mathbf{R}}^n$
be $C_b^{\infty}$.
Then, under Assumptions \textup{(H1)--(H4)}, we have the following
asymptotic expansion as $\varepsilon\searrow0$:
there are real constants $c$ and $\alpha_0, \alpha_1, \ldots$ such that
\begin{eqnarray*}
&&
{\mathbb E} \bigl[
G( Y^{\varepsilon, 1}) \exp\bigl( - F ( Y^{\varepsilon, 1})
/\varepsilon^2 \bigr)
\bigr]
\\
&&\qquad=
\exp\bigl( - F_{\Lambda} (\gamma) /\varepsilon^2 \bigr) \exp(- c/\varepsilon)
\cdot
\bigl(
\alpha_0 +\alpha_1 \varepsilon+ \cdots+ \alpha_m \varepsilon^m
+O(\varepsilon^{m+1})
\bigr)
\end{eqnarray*}
for any $m \ge0$.
\end{tm}

%re2.2 #&#
\begin{re}
The only reason for the boundedness assumption for $\sigma$ and $b$ is
for safety.
It is an important and difficult problem whether Lyons' continuity
theorem holds
for unbounded coefficients under a mild growth condition.
(One of such attempts can be found in~\cite{glj}).
If we have such an extension of the continuity theorem, then Theorem
\ref{thmmain} could
easily be generalized
because localization around~$\gamma$ is crucially used in the proof
(see Section~\ref{secproof} below).
\end{re}

%%%%%%%%%%%%%%%%%%%%%%%%%%%%%%%%%%%%%%%%%%%%%%%%%%%%
%%%%%%%% Tsuika Heuristic Proof
%%%%%%%%%%%%%%%%%%%%%%%%%%%%%%%%%%%%%%%%%%%%%%%%%
%s2.2 #&#
\subsection{A heuristic ``proof''}\label{sec22}
Some readers who are not familiar with Laplace-type asymptotics may
find the argument in this paper too complicated.
So,
in this subsection, we try to help them get a bird's eye view of the
proof of the main theorem.
The argument in this subsection is very heuristic and has no rigorous meaning.\vadjust{\goodbreak}

In this subsection, we denote a generic (rough) path by $w$ (instead of
$w^H$ or $W^H$)
and assume for simplicity that $G \equiv1$ and $\beta$ is independent
of $\varepsilon$ so that $Y^{\varepsilon} = \Psi( \varepsilon w)$
holds at least formally.
As physicists often do,\vspace*{1.5pt}
we will write the law of fBM $\mu^H$ heuristically as
$\mu^H (dw) = Z^{-1} \exp( - |w|^2_{\mathcal{H}^H } /2) \mathcal{D}w$, \vspace*{1pt}where
$\mathcal{D}w$ is the nonexistent ``Lebesgue measure''
and $Z$ is a ``normalizing constant.''

In this case, we study the following quantity:
%
%e2.1 #&#
\begin{eqnarray}\label{heur1}
&&\int\exp\biggl( - \frac{F ( \Psi( \varepsilon w)) }{\varepsilon^2}
\biggr) \mu^H (dw)
\nonumber
\\[-8pt]
\\[-8pt]
\nonumber
&&\qquad=
\frac{1}{Z} \int\exp\biggl( - \frac{ F \circ\Psi( \varepsilon w)
+ |\varepsilon
w|^2_{\mathcal{H}^H}/ 2 }{\varepsilon^2} \biggr) \mathcal{D}w.
\end{eqnarray}
Note that the functional $F_{\Lambda}$ in (H2) appears on the
right-hand side above.
It achieves the minimum at $\gamma$.
So, as in the calculus for freshmen,
one can easily imagine that $(1/2) \times\operatorname{Hessian}$ at $w =\gamma$
plays a very important role.

Let us continue.
By shifting $w \mapsto w +(\gamma/\varepsilon)$, the right-hand side
of (\ref
{heur1}) is equal to
%
%e2.2 #&#
\begin{eqnarray}
&&
\frac{1}{Z} \int\exp\biggl( - \frac{ F \circ\Psi( \varepsilon w +
\gamma)
+ |\varepsilon w + \gamma|^2_{\mathcal{H}^H}/ 2 }{\varepsilon^2} \biggr)
\mathcal{D}w
\nonumber\\
&&\qquad=
\frac{1}{Z} \int\exp\biggl[ - \frac{1}{\varepsilon^2} \biggl\{
F_{\Lambda}
(\gamma)
+ \varepsilon\cdot0 + \frac{\varepsilon^2}{2} ( \langle Aw,w\rangle
+ | w|^2_{\mathcal{H}^H} ) +
O(\varepsilon^3) \biggr\}
\biggr] \mathcal{D}w
\nonumber
\\[-8pt]
\\[-8pt]
\nonumber
&&\qquad=
e^{ - F_{\Lambda} (\Lambda) /\varepsilon^2} \int\exp\biggl[ -
\frac{\langle
Aw,w\rangle}{2} + O(\varepsilon)
\biggr] \mu^H (dw)
\\
&&\qquad\sim e^{ - F_{\Lambda} (\Lambda) /\varepsilon^2} \int\exp
\biggl[ -
\frac
{\langle Aw,w\rangle}{2}
\biggr] \mu^H (dw)\qquad
 \mbox{as $\varepsilon\searrow0$.}\nonumber
\end{eqnarray}
Note that what we did is the Taylor expansion of $F \circ\Psi$ at
$w=\gamma$.
Here, the first order term vanishes, because $\gamma$ is a stationary point.
When the Taylor expansion as above is done, some kind of localization
is usually necessary.
In this case, however, we can localize around $\gamma$, thanks to the
large deviation principle.

As we have seen, the Taylor expansion of $F \circ\Psi$ plays a
central role.
Since we assumed Fr\'echet differentiability of $F$ around $\phi^0
=\Psi(\gamma)$,
the key point is the Taylor expansion of the It\^o map $w \mapsto\Psi
(w)$ around $\gamma$.
This part is rather hard, but was already proved in the author's
previous paper~\cite{ina2}.
See Theorem~\ref{thmtaylor} below, which is a special case of the
result in~\cite{ina2}.

To make this argument rigorous,
we also need
integrability of\break $\exp( - \langle A \bullet, \bullet\rangle/2 )$ with
respect to $\mu^H$.
If $A$ is Hilbert--Schmidt and $A > -\mathrm{Id}$ in the form sense, this
is integrable and its expectation is written
in terms of the Carleman--Fredholm determinant $\det_2 (\mathrm{Id} +A)$.
See Lemma~\ref{lmdet2} and Remark~\ref{recfdet}.
Therefore, it is important to prove the Hilbert--Schmidt property of $A$.
(Precisely, the quadratic form $\langle A \bullet, \bullet\rangle$
must be
replaced by its stochastic extension, i.e.,
the element of the second order Wiener chaos corresponding to the
Hilbert--Schmidt operator $A$.)

%%%%%%%%%%%%%%%%%%%%%%%%%%%%%%%%%%%%%%%%%%%%%%%%%%%%
%%%%%%%% (End) Tsuika Heuristic Proof
%%%%%%%%%%%%%%%%%%%%%%%%%%%%%%%%%%%%%%%%%%%%%
%%%%%%%%%%%%%%%%%%%%%%%%%%%%%%%%%%%%%%%%%%%%%%%%%%%%%%%%%%%%%%%%%%%%%%%
%
%s3 #&#
\section{A review of fractional Brownian rough paths}\label{sec3}
In this section we recall that $d$-dimensional fBM $(w^H_t)_{0 \le t
\le1}$ with Hurst parameter $H \in(1/4, 1/2)$ can be lifted as a
random variable on the geometric rough path space $G\Omega_p ({\mathbf{R}}^d)$ for $1/H < p< [1/H] +1$
(Coutin--Qian~\cite{cq} or Section 4.5 of Lyons--Qian~\cite{lq}).
When $H\in(1/4, 1/3]$, not only the first and the second level paths,
but also the third level paths
play a role.

%
%%%%%%%%%%%%
%s3.1 #&#
\subsection{Geometric rough paths, Lyons' continuity theorem and Taylor expansion of It\^o maps}\label{sec31}
In this subsection we recall definitions of geometric rough paths, and
a rough differential equation (RDE)
and Lyons' continuity theorem for the It\^o map.
We also review (stochastic) Taylor expansion for It\^o maps around a
``nice'' path,
which was shown in~\cite{ina2}.
It plays a crucial role in the proof of the Laplace asymptotic expansion.
Note that no probability measure is involved in this subsection.
No new results are presented in this subsection.

Before introducing the rough path space, let us first introduce some
path spaces in the usual sense and the norms on them.
Let $\mathcal{V}$ be a real Banach space.
Throughout this paper, we assume $\dim\mathcal{V} <\infty$ and the time
interval is $[0,1]$.
In almost all applications in later sections, either $\mathcal{V} ={\mathbf{R}}^d$ or $\mathcal{V} =\operatorname{Mat} (n,d)$
(the space of $n \times d$ matrices).
Let
\[
C = C([0,1 ], \mathcal{V} ) = \{ k \dvtx  [0,1] \to\mathcal{V}  |  \mbox{continuous}\}
\]
be the space of $\mathcal{V}$-valued continuous functions with the usual sup-norm.
%
%For $0< \alpha<1$, let
%$C^{\alpha-hldr} $
%be the set of $k \in C$ such that
%
%Similarly,
%
For $p \ge1$, $C^{p\mbox{-}\mathrm{var}}$ is the set of $k \in C$ such that
\[
\|k\|_{p\mbox{-}\mathrm{var}} ;= |k_0|+ \Biggl(
\sup_{ \mathcal{P}} \sum_{i=1}^n | k_{t_i} - k_{t_{i-1}} |^p \Biggr)^{1/p}
<\infty,
\]
where $\mathcal{P}$ runs over all the finite partition of $[0,1]$.
If $p, q \ge1$ with $1/p + 1/q >1$ and
$k \in C^{q\mbox{-}\mathrm{var}} ( L(\mathcal{V}, \mathcal{W}) )$ and $l \in C^{p\mbox{-}\mathrm{var}} (\mathcal{V})$ with $l_0=0$,
then the Young integral
\[
\int_s^t k_u \,dl_u := \lim_{|\mathcal{P}| \searrow0} \sum_{i=1}^N
k_{t_{i-1}} ( l_{t_{i}} - l_{t_{i-1}} )
\]
is well-defined.
Here, $L(\mathcal{V}, \mathcal{W})$ is the set of linear maps from $\mathcal{V}$
to $\mathcal{W}$
and $\mathcal{P}=\{ s =t_0 <t_1 < \cdots< t_N = t \}$ is a partition of $[s,t]$.
Moreover, $t \mapsto\int_0^t k_u \,dl_u \in{\mathbf{R}}^n$ is of finite
$p$-variation and
$\| \int_0^{\cdot} k_u \,dl_u \|_{p\mbox{-}\mathrm{var}} \le \mathrm{const}\cdot \| k \|_{q\mbox{-}\mathrm{var}} \| l
\|_{p\mbox{-}\mathrm{var}} $.
More precisely, if there is a control function $\omega$ such that
$
| k_t -k_s | \le\omega(s,t)^{1/q},   | l_t -l_s | \le\omega(s,t)^{1/p}$,
then
\[
\biggl| \int_s^t k_u dl_u - k_s (l_t -l_s) \biggr| \le \mathrm{const}\cdot
\omega
(s,t)^{1/p +1/q}.
\]
In particular, if $\tilde{l} \in C^{p\mbox{-}\mathrm{var}} (\mathcal{V})$ and $\tilde{k}
\in C^{q\mbox{-}\mathrm{var}} (\mathcal{W})$
with $1/p +1/q >1$, then
$\int_s^t \tilde{k} _u \otimes d \tilde{l}_u $ is well-defined.

Next we introduce the Besov space $W^{\delta, p}$ for $p > 1$ and
$0<\delta<1$.
For a measurable function $k \dvtx  [0,1] \to\mathcal{V}$, set
%
%e3.1 #&#
\begin{equation} \label{defbesovnorm}
\| k \|_{W^{\delta, p}} = \| k\|_{L^p} +
\biggl(
\iint_{[0,1]^2} \frac{ |k_t - k_s|^p}{ |t-s|^{1 +\delta p}} \,ds\,dt
\biggr)^{1/p}.
\end{equation}
The Besov space $W^{\delta, p}$ is the totality of $k$'s such that $\|
k\|_{W^{\delta, p}} < \infty$.
When $1/p < \delta$, this Banach space is continuously imbedded in $C$
and basically we only consider such a case.
The subspace of functions which start at $0$ (i.e., $k_0 =0$)
is denoted by $C_0$, $C_0^{\alpha\mbox{-}\mathrm{hldr}} $, etc.
When we need to specify the range of functions, we write
$C^{p\mbox{-}\mathrm{var}}(\mathcal{V})$, $W_0^{\delta, p}(\mathcal{V})$, etc.
(The domain is always $[0,1]$ and, hence, is usually omitted.)

%
%%%%%%%%%%%%%

Now we introduce the geometric rough path space.
Let $p \ge1$ for a while. (In later sections, however, only the case
$2<p<4$ will be considered.)
Set $\triangle=\{ (s,t)  |  0 \le s \le t \le1\}$.
The $p$-variation norm of a continuous map $A$ form $\triangle$ to a real
finite dimensional Banach space $\mathcal{V}$ is
defined by
\[
\|A\|_{p\mbox{-}\mathrm{var}} = \Biggl( \sup_{ \mathcal{P}} \sum_{i=1}^n | A_{t_{i-1} ,
t_{i}} |^p_\mathcal{V} \Biggr)^{1/p},
\]
where $\mathcal{P}$ runs over all the finite partition of $[0,1]$.
A continuous map
\[
X =\bigl(1, X^1, X^2, \ldots, X^{[p]}\bigr) \dvtx  \triangle\to T^{[p]} (\mathcal{V})=
{\mathbf{R}} \oplus\mathcal{V} \oplus\mathcal{V}^{\otimes2} \cdots\oplus
\mathcal{V}^{\otimes[p]}
\]
is said to be a $\mathcal{V}$-valued rough path of roughness $p$ if it
satisfies the following conditions:

\begin{longlist}[(a):]
\item[(a):] For any $s \le u \le t$, $X_{s,t} = X_{s,u} \otimes X_{u,t}
$, where $\otimes$ denotes the tensor
operation in the truncated tensor algebra $T^{[p]} (\mathcal{V})$.\vspace*{-1pt}
In other words,
$X_{s,t}^j = \sum_{ i=0}^j X_{s,u}^i \otimes X_{u,t}^{j-i}$ for all $1
\le j \le[p]$.
This is called Chen's identity.

\item[(b):]
For all $1 \le j \le[p]$, $\| X^j \|_{p/j \mbox{-}\mathrm{var}} <\infty$.
\end{longlist}

We usually omit the $0$th component $1$ and simply write $X=( X^1,
\ldots, X^{[p]}) $.
The first level path of $X$ is naturally regarded as an element in
$C^{p\mbox{-}\mathrm{var}}_0 (\mathcal{V})$
by $t \mapsto X^1_{0,t}$.\vspace*{1pt} [We will abuse the notation to write $X^1 \in
C^{p\mbox{-}\mathrm{var}}_0 (\mathcal{V})$, e.g.]
The set of all the $\mathcal{V}$-valued rough paths of roughness $p$ is
denoted by $\Omega_p (\mathcal{V})$.
With the distance $d_p (X,Y) = \sum_{i=1}^{[p]} \| X^j -Y^j \|_{p/j
\mbox{-}\mathrm{var}} $,
it becomes a complete metric space.

A $\mathcal{V}$-valued finite variational path $x \in C^{1\mbox{-}\mathrm{var}}_0 (\mathcal{V})$ is naturally lifted
as an element of $\Omega_p (\mathcal{V})$ by the following iterated
Stieltjes integral:
%
%e3.2 #&#
\begin{equation}\label{eqiterint}
X^j_{s,t} = \int_{s \le t_1 \le\cdots\le t_j \le t} \,dx_{t_1}
\otimes dx_{t_2} \otimes\cdots\otimes dx_{t_j}.
\end{equation}
We say $X$ is the smooth rough path lying above $x$.
It is well known that the injection $x \mapsto X \in\Omega_p (\mathcal{V})$ is continuous with respect to the 1-variation norm.
The space of geometric rough path $G\Omega_p (\mathcal{V})$ is the closure
of $C_0^{1\mbox{-}\mathrm{var}} (\mathcal{V})$ with respect to $d_p$.
Since $\mathcal{V}$ is separable,
$G\Omega_p (\mathcal{V})$ is a complete separable metric space.

%%%%%%%%%%%%%%%%%%%%%%%%%%%%%%%%%
Let us recall some properties of the $q$-variational path for $1 \le q <2$.
For the facts presented below, see Section 3.3.2 in Lyons--Qian \cite
{lq} or Inahama~\cite{ina2}, for example.
Since $k \in C_0^{q\mbox{-}\mathrm{var}} (\mathcal{V})$ is Young integrable with respect
to itself,
the iterated integral in (\ref{eqiterint}) is still well-defined and
$k$ can be lifted to an element $K \in G\Omega_p (\mathcal{V})$ if $p \ge2$.
This injection $C_0^{q\mbox{-}\mathrm{var}} (\mathcal{V}) \hookrightarrow G\Omega_p
(\mathcal{V})$ is continuous.

For any $m=1,2,\ldots$ and any $k \in C_0 (\mathcal{V})$, the $m$th dyadic
piecewise linear approximation $k(m)$ is defined by
\[
k(m)_t = k_{(l-1) /2^m} + 2^m \bigl( k_{l /2^m} -k_{(l-1) /2^m} \bigr) \bigl(t-
(l-1)/2^m \bigr)
\qquad
\mbox{if }t \in\biggl[ \frac{l-1}{2^m}, \frac{l}{2^m} \biggr].
\]
If $k$ is of $q$-variation ($q \ge1$), then $k(m)$ converges to $k$ in
$(q +\varepsilon)$-variation norm
for any $\varepsilon>0$.
It implies that, if $p \ge2$ and $k \in C_0^{q\mbox{-}\mathrm{var}} (\mathcal{V})$ for $1
\le q <2$, then $K(m)$ converges to $K$
in $G\Omega_p (\mathcal{V})$.

Suppose that if $p \ge2, 1 \le q <2$, and $1/p +1/q >1$, then the shift
\[
(X, k) \in G\Omega_p (\mathcal{V}) \times C_0^{q\mbox{-}\mathrm{var}} (\mathcal{V}) \mapsto
X + K \in G\Omega_p (\mathcal{V})
\]
is well-defined by the Young integral and this map is continuous.
Similarly,
\[
(X, k) \in G\Omega_p (\mathcal{V}) \times C_0^{q\mbox{-}\mathrm{var}} (\mathcal{W}) \mapsto
(X , K) \in G\Omega_p (\mathcal{V} \oplus\mathcal{W})
\]
is well-defined and continuous.
These facts are well known.
(See Section 9-4 in Friz and Victoir's book~\cite{fvbook}, e.g.)
Note that the notation ``$X+K$'' above may be somewhat misleading since
the geometric rough path space is not an additive group.

%%%%%%%%%%%%%%%%%%%%%%%%%%%%%%%%%%%%%%%%%%%%%%%%%%
%%%%%%%%%%%%%%%%%%%%%%%%%%%%%%%%%%%%%%%%%%%%%%%%%%%%

Let $\mathcal{V}$ and $\mathcal{W}$ be two finite dimensional real Banach spaces
and let $\sigma\dvtx  \mathcal{W} \to L(\mathcal{V}, \mathcal{W} )$ with some
regularity condition, which will be specified later.
We consider the following differential equation in the rough path sense
(rough differential equation or RDE):
%
%e3.3 #&#
\begin{equation}\label{defrdekihon}
dY_t = \sigma(Y_t ) \,dX_t ,\qquad   Y_0=y_0 \in\mathcal{W}.
\end{equation}
When there is a unique solution $Y$ for given $X$,
it is denoted by $Y=\Phi(X)$ and the map $\Phi\dvtx  G\Omega_p (\mathcal{V})
\to G\Omega_p (\mathcal{W})$ is called the It\^o map.\vadjust{\goodbreak}

The following is called Lyons' continuity theorem (or universal limit theorem)
and is most important in the rough path theory.
(See Section 6.3, Lyons--Qian~\cite{lq}.
For a proof of continuity when the coefficient $\sigma$ also varies,
see Inahama~\cite{ina2}, e.g.)
%
%th3.1 #&#
\begin{tm}
\label{thmconti}
\textup{(i)}
Let $p \ge2$ and assume that $\sigma\in C^{[p]+1}_b ( \mathcal{V}, \mathcal{W} )$.
Then, for given $X \in G\Omega_p (\mathcal{V})$ and an initial value $y_0
\in\mathcal{W}$, there is a unique solution $Y \in G\Omega_p (\mathcal{W})$
of RDE (\ref{defrdekihon}).
Moreover, there is a constant $C_M>0$ for $M>0$ such that,
if
\[
|y_0| \le M, \qquad \sum_{j=1}^{[p]} \|X^j\|_{p/j \mbox{-}\mathrm{var}} \le M,\qquad
 \sum_{j=0}^{[p]+1} \sup_{y \in\mathcal{W}} \| \nabla^j \sigma(y)
\| \le M,
\]
then $\sum_{j=1}^{[p]} \|Y^j\|_{p/j \mbox{-}\mathrm{var}} \le C_M$.\vspace*{-4pt}

\begin{longlist}[(ii)]
\item[(ii)]
Keep the same assumption as above.
Assume that $X_l \to X$ in $G\Omega_p (\mathcal{V}) $ and $y_0^l \to y_0$
in $\mathcal{W}$ as $l \to\infty$.
Assume further that $\sigma_l, \sigma\in C^{[p]+1}_b ( \mathcal{V},
\mathcal{W} )$ satisfy that
\[
\sup_{l \ge1} \sum_{j=0}^{[p]+1} \sup_{y \in\mathcal{W}} \| \nabla^j
\sigma_l (y) \| \le M
\]
for some constant $M>0$ and
\[
\lim_{l \to\infty}
\sum_{j=0}^{[p]+1} \sup_{|y|_{\mathcal{W}} \le N } \| \nabla^j \sigma_l
(y) -\nabla^j \sigma(y) \| =0
\]
for each fixed $N>0$.
Then, $Y_l \to Y$ in $G\Omega_p (\mathcal{W}) $,
where $Y_l$ is the solution of RDE~(\ref{defrdekihon}) corresponding
to $(X_l, y_{0}^{l} , \sigma_l)$.
\end{longlist}
\end{tm}

%%%%%%%%%%%%%%%%%%%%%%%%%%
%TAYLOR TAYLOR !!!!!
%%%%%%%%%%%%%%%%%%%%%%%%%
%
In this paper we consider the following RDE indexed by small parameter
$\varepsilon>0$.
Let $\sigma\in C_b^{\infty} ( {\mathbf{R}}^n , \operatorname{Mat}(n,d))$ and $\beta
\in C_b^{\infty} ([0,1] \times{\mathbf{R}}^n, {\mathbf{R}}^n)$.
For fixed $\varepsilon\in[0,1]$, consider
%
%e3.4 #&#
\begin{equation}
\label{defrde}
dY^{\varepsilon}_t = \sigma(Y^{\varepsilon}_t ) \varepsilon \,dX_t
+\beta(\varepsilon, Y^{\varepsilon}_t) \,dt,\qquad
  Y^{\varepsilon}_0=0.
\end{equation}
[This is the same RDE as in (\ref{defrdeintro})].
If we define $\hat\sigma_{\varepsilon}\dvtx  {\mathbf{R}}^n , \operatorname{Mat}(n,d+1)$ by
\[
\hat\sigma_{\varepsilon} (y) x' = \sigma(y) x + \beta(\varepsilon
, y) x_{d+1},
\qquad
x'=(x ,x_{d+1}) \in{\mathbf{R}}^d \oplus{\mathbf{R}},
\]
then $Y^{\varepsilon} = \hat\Phi_{\varepsilon} (\varepsilon X,
\lambda)$.
Here, $\lambda_t =t$ and $ \hat\Phi_{\varepsilon}\dvtx  G\Omega_p ({\mathbf{R}}^{d+1})
\to G\Omega_p ({\mathbf{R}}^n )$ is the It\^o map which corresponds to
$\hat
\sigma_{\varepsilon}$.
Note that $\hat\sigma_{\varepsilon}$ converges to $\hat\sigma
_{\varepsilon'}$ in the
sense of Theorem~\ref{thmconti}(ii)
as $\varepsilon\to\varepsilon'$.

Now we consider the (stochastic) Taylor expansion around $\gamma\in
C_0^{q\mbox{-}\mathrm{var}} ({\mathbf{R}}^d)$
with $1/p+1/q >1$.
Consider $\hat\Phi_{\varepsilon} (\varepsilon X +\gamma, \lambda)$
or, equivalently,
the solution of
the following RDE:
%
%e3.5 #&#
\begin{equation}\label{eqshiftedrde}
d \tilde{Y}^{\varepsilon}_t = \sigma(\tilde{Y}^{\varepsilon}_t )
(\varepsilon \,dX_t +
d\gamma
_t) +\beta(\varepsilon, \tilde{Y}^{\varepsilon}_t) \,dt,\qquad
\tilde{Y}^{\varepsilon}_0=0.\vadjust{\goodbreak}
\end{equation}
We will write $\phi^{(\varepsilon)} = ( \tilde{Y}^{\varepsilon}
)^1$ (the first
level path).
Note that $\hat\Phi_{0} (\gamma, \lambda)$ is lying above $\phi^0
= \Psi
(\gamma) \in C_0^{q\mbox{-}\mathrm{var}} ({\mathbf{R}}^n)$
which is defined by
%
%e3.6 #&#
\begin{equation}\label{eqrdephi0}
d \phi^0_t = \sigma(\phi^0_t ) \,d\gamma_t +\beta(0 , \phi^0_t) \,dt,\qquad
 \phi^{0}_0=0.
\end{equation}

In the following theorem, we consider the asymptotic expansion of
$\phi^{(\varepsilon)}- \phi^0$.
By formally operating $(m!)^{-1} (d/d\varepsilon)^m \vert
_{\varepsilon=0}$ on both
sides of (\ref{eqshiftedrde}),
we get an RDE for the $m$th term $\phi^m$ (see~\cite{ina2} for detail).
Note that $\phi^m$ depends on $X, \gamma$ (although~$\gamma$ is
basically fixed in this paper), but independent of $\varepsilon$.
(The superscript $m$ does not denote the level of the path $\phi^m$.
Here we only consider the usual paths or the first level paths.)

In what follows, we will use the following notation;
for a geometric rough path~$X$ of roughness $p$,
%
%e3.7 #&#
\begin{equation}
\label{defxi}
\xi(X)= \| X^1 \|_{p\mbox{-}\mathrm{var}} + \| X^2 \|_{p/2\mbox{-}\mathrm{var}}^{1/2} +\cdots+ \bigl\|
X^{[p]} \bigr\|_{p/[p]\mbox{-}\mathrm{var}}^{1/[p]}.
\end{equation}

%%%%%%%%%%%%%%%%%%%%%%%%%%%%%%%%%%%
%
%th3.2 #&#
\begin{tm}
\label{thmtaylor}
Let $p \ge2$, $1 \le q <2 $ with $1/p+1/q>1$ and let the notation be
as above.
Then, for any $m =1,2,\ldots,$
we have the following expansion:
\[
\phi^{(\varepsilon)} = \phi^0 +\varepsilon\phi^1 + \cdots+
\varepsilon^m \phi^m +
R^{m+1}_{\varepsilon}.
\]
The maps $(X, \gamma) \in G\Omega_p ({\mathbf{R}}^d) \times C_0^{q\mbox{-}\mathrm{var}}
({\mathbf{R}}^d) \mapsto\phi^k, R^{m+1}_{\varepsilon} \in C_0^{p\mbox{-}\mathrm{var}}
({\mathbf{R}}^n) $
are continuous $(0 \le k \le m)$.
Moreover, the following estimates \textup{(a)}, \textup{(b)} hold:

\begin{longlist}[(a)]
\item[(a)]  For any $r_1>0$, there exists $C_1 >0$ which depends
only on $r_1$
such that,
if $\|\gamma\|_{q\mbox{-}\mathrm{var}} \le r_1$, then $ \|\phi^k \|_{p\mbox{-}\mathrm{var}} \le C_1 (1+
\xi(X))^k$ holds.

\item[(b)]
For any $r_2, r_3>0$, there exists $C_2 >0$ which depends only on
$r_2, r_3$
such that,
if $\|\gamma\|_{q\mbox{-}\mathrm{var}} \le r_2$ and $\xi(\varepsilon X) \le r_3$,
then $ \| R_{\varepsilon}^{m+1} \|_{p\mbox{-}\mathrm{var}} \le C_2 (\varepsilon+\break \xi
( \varepsilon X))^{m+1}$ holds.
\end{longlist}
\end{tm}

%%%%%%%%%%%%%%%%%%%%%%%%%%%%%%%%%%%%%%%
%s3.2 #&#
\subsection{Fractional Brownian rough paths}\label{sec32}

First
we introduce fractional Brownian motion (fBM for short) of Hurst
parameter $H$.
There are several books and surveys on fBM (see~\cite{bhoz,cou,mis},
e.g.).
In this paper we only consider the case $1/4 <H < 1/2$.
A real-valued continuous stochastic process $(w^H_t)_{t \ge0}$
starting at $0$ is
said to a fBM of Hurst parameter $H$ if it is a centered Gaussian process
with
\[
{\mathbb E} [ w^H_t w^H_s ] =\tfrac12 [ t^{2H} +s^{2H} -|t-s|^{2H}]
\qquad
(s,t \ge0).
\]
This process has stationary increments
$
{\mathbb E} [ (w^H_t -w^H_s )^2 ] =|t-s|^{2H}
\ (s,t \ge0),\vspace*{1pt}
$
and self-similarity, that is, for any $c>0$,
$
( c^{-H} w^H_{ct} )_{t \ge0}
$
and
$( w^H_{t} )_{t \ge0} $ have the same law.
Note that $(w^{1/2}_t)_{t \ge0}$ is the standard Brownian motion.
For $d \ge1$,
a $d$-dimensional fBM\vspace*{1pt} is defined by $( w^{H,1}_{t} , \ldots,w^{H,d}_{t}
)_{t \ge0}$,
where $w^{H,i}\ (i=1,\ldots,d)$ are independent one-dimensional fBM's.\vspace*{1pt}
Its law $\mu^H$ is a probability measure on $C_0({\mathbf{R}}^d)$.
[Actually, it is a probability measure on $C_0^{p\mbox{-}\mathrm{var}} ({\mathbf{R}}^d)$ for
$p >1/H$.]

%%%%%%%%%%%%%%%%%%%%%%%%%%%%%%%%%%%

%Let $K^H$ be the Volterra kernel given by
%K^H (t,s) = \frac{ (t-s)^{H-1/2}}{\Gamma(H + 1/2)} F (
%H-\frac12, \frac12 -H, H+\frac12 , 1 -\frac{s}{t}
%)
%{\bf1}_{\{s<t\}}.
%Here, $F$ denotes the Gauss hypergeometric function.
%For a $d$-dimensional standard Brownian motion $(b_t)_{t \ge0}$, the
%process
%t \mapsto\int_0^t K^H(t,s) \,db_s
%becomes a $d$-dimensional fBm with Hurst parameter $H$ (see
%Decreusefond-\"Ustsunel~\cite{du}).

%%%%%%%%%%%%%%%%%%%%%%%%%%%%%%%%%%%%%%%%
Let $H \in(1/4,1/2)$.
We denote by $w^H(m)$ the $m$th dyadic piecewise linear approximation
of $w^H$, that is,
piecewise linear approximation associated with the partition $\{ j
2^{-m}  |  0 \le j \le2^m \}$.
The existence of a fractional Brownian rough path (fBRP for short) was
shown by Coutin--Qian~\cite{cq}
as an almost sure limit of $W^H (m)$ as $m \to\infty$,
where $W^H(m)$ is the smooth rough path lying above $w^H(m) \in
C_0^{1\mbox{-}\mathrm{var}} ({\mathbf{R}}^d)$.
More precisely, they proved
\[
{\mathbb E} \Biggl[
\sum_{m=1}^{\infty} \| W^H(m+1)^j - W^H(m)^j \|_{p/j \mbox{-}\mathrm{var}}
\Biggr] < \infty
\qquad
(1 \le j \le[p]).
\]
In particular, $W^H(m)$ converges to $W^H$ in the $L^1$-sense, too.
When $1/3 <H< 1/2$, $[p]=2$ and when $1/4 <H \le1/3$, $[p]=3$.

Now we prove a theorem of Fernique-type for fBRP for later use.
We give a direct proof here for readers' convenience by
using a useful estimate in Millet and Sanz-Sole~\cite{mss}.
(The case $H=1/2$ is shown in~\cite{ina}, e.g.)
It should be noted, however, that (i) this proposition is included in
Theorems~15.22 and 15.42,~\cite{fvbook} and
(ii) Friz and Oberhauser~\cite{fo} recently showed this kind of integrability
for a wider class of Gaussian rough paths, by using isoperimetric inequality.
%
%pr3.3 #&#
\begin{pr}
\label{prfern}
Let $1/4<H <1/2$ and $W^H$ be a $d$-dimensional fBRP as above.
\begin{longlist}[(1)]
\item[(1)]
Then, there exists a positive constant $c$ such that
\[
{\mathbb E} [ \exp( c \xi( W^H)^2 ) ] = \int_{ G\Omega_p ({\mathbf{R}}^d)}
\exp( c \xi( X)^2 ) {\mathbb P}^H (dX) <\infty,
\]
where $\xi$ is given in (\ref{defxi}) and ${\mathbb P}^H $ denotes the
law of $W^H$.

\item[(2)] For any $r >0$ and $1 \le j \le[p]$,
\mbox{$\lim_{m \to\infty} {\mathbb E} [ \| W^H(m)^j - W^{H,j} \|^r_{p/j
\mbox{-}\mathrm{var}} ] = 0$}.
\end{longlist}
\end{pr}

\begin{pf}
In this proof, $c_1, c_2,\ldots$ are positive constants which may
change from line to line.
For a rough path $X$ of roughness $p$ and $\gamma>p-1$, set
\begin{eqnarray}
D_{j,p} (X, Y) = \Biggl(
\sum_{n=1}^{\infty} n^{\gamma} \sum_{l=1}^{2^n} \bigl|X^j_{(l-1)/2^n, l/2^n}
- Y^j_{(l-1)/2^n, l/2^n} \bigr|^{p/j}
\Biggr)^{j/p}
\nonumber\\
\eqntext{(1 \le j \le[p]).}
\end{eqnarray}
When $Y=0$, we write $D_{j,p} (X) =D_{j,p} (X, Y) $ for simplicity.
From Section 4.1 in Lyons--Qian~\cite{lq}, the following estimates hold:
%
%e3.8 #&#
\begin{eqnarray}\label{ineqDjp}
&&\| X^1 -Y^1 \|_{p \mbox{-}\mathrm{var}}^{p} \le c_1 D_{1,p} (X, Y)^p
\nonumber\\
&&\| X^2 -Y^2 \|_{p/2 \mbox{-}\mathrm{var}} ^{p/2}\nonumber\\
 &&\qquad\le c_1
\bigl[
D_{2,p} (X, Y)^{p/2} + D_{1,p} (X, Y)^{p/2} \bigl( D_{1,p} (X)^p+D_{1,p} (Y)^p\bigr)^{1/2}
\bigr]
\nonumber\\
\qquad &&\| X^3 -Y^3 \|_{p/3 \mbox{-}\mathrm{var}}^{p/3} \\
&&\qquad\le c_1
\bigl[
D_{3,p} (X, Y)^{p/3} + D_{2,p} (X, Y)^{p/3} \bigl( D_{1,p} (X)^p+D_{1,p} (Y)^p\bigr)^{1/3}
\nonumber\\
&&\qquad\qquad{}
+ D_{1,p} (X, Y)^{p/3} \bigl( D_{2,p} (X)^{p/2} + D_{2,p} (Y)^{p/2}\bigr)^{2/3}
\nonumber\\
&&\hspace*{78pt}\qquad\quad{}
+ D_{1,p} (X, Y)^{p/3} \bigl( D_{1,p} (X)^{p} + D_{1,p} (Y)^{p}\bigr)^{2/3}
\bigr].\nonumber
\end{eqnarray}

Proposition 2 in~\cite{mss} states that there is a sequence $\{ a_m\}$
of positive numbers converging to $0$
such that, for any $r >p$,
\[
{\mathbb E} [ D_{j,p} (W^H(m) , W^H)^{r} ]^{1/r} \le a_m r^{j/2}
\]
holds.
For simplicity, set $F_m =D_{j,p} (W^H(m) , W^H)^{2/j}$.
Then, from the above inequality,
\[
{\mathbb P} (N< F_m ) \le N^{-N} {\mathbb E} [F_m^N] \le c_2^N a_m^N
\]
for $N=4,5,\ldots.$
Therefore,
\begin{eqnarray*}
{\mathbb E} [e^{cF_m} ]
&\le&
\sum_{N=0}^{\infty} e^{c(N+1)} {\mathbb P} (N< F_m \le N+1)
\\
&\le&
(e^c +\cdots+c^{4c})
+
e^c \sum_{N=4}^{\infty} e^{cN} {\mathbb P} (N< F_m )
\\
&\le&
(e^c +\cdots+c^{4c})
+
e^c \sum_{N=4}^{\infty} \exp[ N (c+ \log c_2 - \log a_m) ].
\end{eqnarray*}
For given $c>0$, there exists $m_0$ such that $m \ge m_0$ implies $c+
\log c_2 - \log a_m <0$.
Thus, we obtain
\[
\sup_{m \ge m_0} {\mathbb E} [e^{cF_m} ]
\le
\sup_{m \ge m_0} {\mathbb E} [ \exp(c D_{j,p} (W^H(m) , W^H)^{2/j} )]
<\infty.
\]
On the other hand, it is easy to see that, for each fixed $m_0$,
there is a constant $c'(m_0) >0$ such that
$D_{j,p} (W^H(m_0) )^{1/j} \le c'(m_0) \|w^H\|_{\infty}$.
Hence, the usual Fernique theorem for Gaussian measures applies and
$D_{j,p} (W^H(m_0) )^{1/j} $
is square exponentially integrable.
Using (\ref{ineqDjp}) and the triangle inequality for $D_{j,p}$, we
prove~(1).
In a similar way, we see that
\[
\sup_{m \ge1} {\mathbb E} [ D_{j,p} (W^H(m) )^{r} ] <\infty,\qquad
\sup_{m \ge1} {\mathbb E} [ \| W^H(m)^j \|^{r}_{p/j \mbox{-}\mathrm{var}} ]
<\infty.
\]
This implies (2).
\end{pf}

%%%%%%%%%%%%%%%%%%%%%%%%%%%%%%%%%%%%%%%%%%%%%%%%
%
%%%%%%%%%%%%%%%%%%%%%%%%%%%%%%%%%%%%%%

Let $\mathcal{H}^H$ be the Cameron--Martin subspace of fBM
[i.e., $k \in C_0 ({\mathbf{R}}^d)$ is an element of $\mathcal{H}^H$ if and
only if
$\mu^H$ and $\mu^H(  \cdot  +k)$ are mutually absolutely continuous].\vadjust{\goodbreak}
%
%For $h \in L^2 =L^2 ([0,1], {\mathbf{R}}^d)$, set
%It is known that
%$k: =Uh \in\mathcal{H}^H$ and the map $h \mapsto k=Uh$ is unitary from
%$L^2$ to $\mathcal{H}^H$.
%
When $H=1/2$, it is easy to see $k \in\mathcal{H}^{1/2}$ is of finite
$1$-variation.
But, when $H \in(1/4,1,2)$, does $k \in\mathcal{H}^{H}$ have a similar
nice property in terms of variation norm?
The following theorem answers this question.
As a result, $\mathcal{H}^{H}$ is continuously (and compactly) embedded in
$G\Omega_p ({\mathbf{R}}^d)$ for $p \ge2$.
%
%pr3.4 #&#
\begin{pr}[(Friz--Victoir~\cite{fv})]\label{prFV}
\begin{longlist}[(ii)]
\item[(i)]
Let $0<\delta<1$ and $p \ge1$ such that $\alpha= \delta-1/p >0$ and
set $q= 1/\delta$.
Then, we have a continuous embedding
\[
W^{\delta, p} \subset C^{q\mbox{-}\mathrm{var}} ,  \qquad W^{\delta, p} \subset
C^{\alpha\mbox{-}\mathrm{hldr}}.
\]
More precisely, for $h \in W^{\delta, p} $,
\[
\omega(s,t) = \|h\|^q_{W^{\delta, p} ; [s,t]} (t-s)^{\alpha q},
\qquad
0 \le s \le t \le1
\]
becomes a control function in the sense of Lyons--Qian~\cite{lq}, page
16, and $h$ is controlled by a constant multiple of $\omega$
[i.e., $|h_t -h_s| \le \mathrm{const} \times\omega(s,t)^{1/q}$].

\item[(ii)] Let the Hurst parameter $H \in(0, 1/2)$.
If $1/2 <\delta< H+ 1/2$, then $\mathcal{H}^H \Subset W_0^{\delta, 2}$
(compact embedding).
Therefore, for any $\alpha\in(0,H)$ and $q \in((H+1/2)^{-1} , 2)$,
\[
\mathcal{H}^H \Subset C_0^{\alpha\mbox{-}\mathrm{hldr}},\qquad
\mathcal{H}^H \Subset C_0^{q\mbox{-}\mathrm{var}}.
\]
\end{longlist}
\end{pr}

We give a theorem of a Cameron--Martin type for fBRP $W^H$.
(For BRP, see~\cite{ina}, e.g.)
Let $1/4 <H< 1/2$ and $1/H <p < [1/H] +1$. Then, fBRP $W^H$ exists on
$G\Omega_p ({\mathbf{R}}^d)$
and its law is a probability measure on $G\Omega_p ({\mathbf{R}}^d)$.
By Proposition~\ref{prFV}, there exists $1 \le q<2$ such that $ \mathcal{H}^H \Subset C_0^{q\mbox{-}\mathrm{var}} \subset G\Omega_p ({\mathbf{R}}^d)$ and $1/p +1/q >1$.
Hence, the shift $X\mapsto X +K $ for $k \in\mathcal{H}^H$ is
well-defined in $G\Omega_p ({\mathbf{R}}^d)$,
where $K$ is the lift of $k$ as usual.
%
%pr3.5 #&#
\begin{pr}
\label{prCM}
Let $\varepsilon>0$ and let ${\mathbb P}_{\varepsilon}^H$ be the law
of $\varepsilon W^H$.
Then, for any $k \in\mathcal{H}^H$,
${\mathbb P}_{\varepsilon}^H$ and ${\mathbb P}_{\varepsilon}^H (
\cdot  +K)$ are
mutually absolutely continuous
and, for any bounded Borel function $f$ on $G\Omega_p ({\mathbf{R}}^d)$,
\begin{eqnarray*}
&&\int_{G\Omega_p ({\mathbf{R}}^d)} f(X+K) {\mathbb P}_{\varepsilon}^H (dX)
\\
&&\qquad=
\int_{G\Omega_p ({\mathbf{R}}^d)} f(X)
\exp\biggl(
\frac{1}{\varepsilon^2} \langle k , X^1 \rangle-\frac
{1}{2\varepsilon^2} \|k\|^2_{\mathcal{H}^H}
\biggr)
{\mathbb P}_{\varepsilon}^H (dX).
\end{eqnarray*}
Here, $ \langle k , X^1 \rangle$ is the measurable linear functional
associated with $k \in\mathcal{H}^H=(\mathcal{H}^H)^*$
for the fBM $t \mapsto X^1_{0,t}$ (i.e., the element of the first
Wiener chaos of the fBM $ X^1$ associated with $k$).
\end{pr}

\begin{pf}
Since $W^H(m) \to W^H$ in $G\Omega_p ({\mathbf{R}}^d)$ and
$k(m) \to k$ in $q$-variation norm as $m \to\infty$, respectively,
$W^H +K = \lim_{m\to\infty} [ W^H(m) + K(m)]$.
On the other hand, $W^H(m) + K(m)$ is the lift of $w^H (m)+k (m)=(w^H + k)(m)$.
Hence, the problem reduces to the usual Cameron--Martin theorem for fBM~$w^H$.\vadjust{\goodbreak}
\end{pf}

In the end of this subsection we give a Schilder-type large deviation
principle for the law of $\varepsilon W^H$
as $\varepsilon\searrow0$.
This was shown by Millet and Sanz-Sole~\cite{mss} (and by Friz--Victoir
\cite{fv,fv0}).

%pr3.6 #&#
\begin{pr}\label{prldp}
Let ${\mathbb P}_{\varepsilon}^H$ be the law of $\varepsilon W^H$ as
above $(1/4 < H <1/2)$.
As before, $1/H <p< [1/H] +1$ and $G\Omega_p ({\mathbf{R}}^d)$ is equipped
with the $p$-variation metric.
Then, as $\varepsilon\searrow0$, $\{ {\mathbb P}_{\varepsilon}^H \}
_{\varepsilon>0}$
satisfies a large deviation principle
with a good rate function $I$, which is given by
\[
I(X) =
\cases{
\frac12
\| k \|_{\mathcal{H}^H}^2 & \quad $\mbox{(if $X$ is lying above $k \in\mathcal{H}^H$),}$ \vspace*{2pt}\cr
\infty&\quad $\mbox{(otherwise).}$ }
\]
\end{pr}

%%%%%%%%%%%%%%%%%% Hilbert--Schmidt property of Hessian
%%%%%%%%%%%%%%%%%%%%%%%%%%%%%%%%
%%%%%%%%%%%%%%%%%%%%%%%%%%%%%%%%%%%%%%%%%%%%%%%%%%%%%%%%%%%%%%%%%%%%%%
%%%
%s4 #&#
\section{Hilbert--Schmidt property of Hessian}\label{sec4}
%%%%%%%%%%%%%%%%%%%%%

In this section we consider the It\^o map restricted on the
Cameron--Martin space $\mathcal{H}^H$ of the
fBM with Hurst parameter $H \in(1/4, 1/2)$ and prove that its Hessian
is a symmetric Hilbert--Schmidt bilinear form.

Throughout this section
we set
%$\sigma(y)=\sigma(y)$ and
$\beta_0 (y)=\beta(0,y)$ for simplicity.
Consider the following RDE:
%
%e4.1 #&#
\begin{equation} \label{rdeqvar}
dY_t = \sigma(Y_t) \,dX_t + \beta_0 (Y_t) \,dt,\qquad   Y_0=0.
\end{equation}
The It\^o map $X \in G\Omega_p( {\mathbf{R}}^d) \mapsto\hat\Phi_0 (X,
\lambda) =Y \in G\Omega_p ( {\mathbf{R}}^n) $
restricted on the Cameron--Martin space $\mathcal{H}^H$ of fBM is denoted
by $\Psi$, that is,
$\Psi(k) =\hat\Phi_0 (K, \lambda)$ for $k \in\mathcal{H}^H$.
Here, $K$ is a geometric rough path lying above $k$ and $\lambda_t =t$.
(Since $k$ is of finite $q$-variation for some $q <2$, as we will see
below, this is well-defined.
Regularity of $k \in\mathcal{H}^H$ in a $p$-variational setting is
studied by Friz--Victoir~\cite{fv}.
Fortunately, $h$ is of finite $q$-variation for some $q<2$ and, hence,
the Young integral is possible.)

%%%%%%%%%%%%% %%%%%%%%%%%%%%%%%%%%%%%%%%

The aim of this section is to prove the following theorem.
Let $F$ and $p'$ be as in Assumption (H1).
%
%th4.1 #&#
\begin{tm}\label{thmHSmain}
$\nabla^2 (F \circ\Psi)(\gamma) \langle \cdot  ,  \cdot
\rangle$
is a
symmetric Hilbert--Schmidt bilinear form on $\mathcal{H}^H$
for any $\gamma\in\mathcal{H}^H$.
\end{tm}

%%%%%%%%%%%%%%%%%%%%%%%%%%%%%%%%%%%%%%%

%re4.2 #&#
\begin{re}
The reader may find arguments in this section a little bit messy.
So, we give a brief summary here.
The most difficult part in proving the above theorem is to show that
the bilinear functional
\[
(f,k) \in\mathcal{H}^H \times\mathcal{H}^H \mapsto\int_0^{\cdot} f_u
\otimes dk_u \in C_0^{p \mbox{-}\mathrm{var}} ({\mathbf{R}}^d \otimes{\mathbf{R}}^d)
\]
is ``Hilbert--Schmidt.''
To compute the Hilbert--Schmidt norm, we need a simple orthonormal basis.
However, we do not know a good basis of $\mathcal{H}^H$.
Therefore, we first imbed $\mathcal{H}^H$ into a larger Hilbert space
$L^{\delta, 2}_{\mathrm{real}}$,
since it has a very simple orthonormal basis of cosine functions, and
then prove the Hilbert--Schmidt property for the norm of $L^{\delta, 2}_{\mathrm{real}}$.
[See (\ref{ineqparapq}) below for the definition of $p$ and $\delta=1/q$.]
\end{re}

%%%%%%%%%%%%%%%%%%%%%%

Note that
\begin{eqnarray*}
\nabla^2 (F \circ\Psi)(\gamma) \langle f,k \rangle
&=&
\nabla F ( \Psi(\gamma) ) \langle\nabla^2 \Psi(\gamma)
\langle f ,k
\rangle\rangle
\\
&&{}
+
\nabla^2 F ( \Psi(\gamma) ) \langle\nabla\Psi(\gamma)
\langle f
\rangle
, \nabla\Psi(\gamma) \langle k \rangle\rangle.
\end{eqnarray*}
ODEs for $\nabla\Psi(\gamma) \langle k \rangle$ and $ \nabla^2
\Psi
(\gamma)
\langle f ,k \rangle$
will be given in (\ref{defchik})--(\ref{defpsi}) below.

Now we set conditions on parameters.
First we have the Hurst parameter $H \in(1/4, 1/2)$.
Then, we can choose $p$ and $q=\delta^{-1}$ such that
%
%e4.2 #&#
\begin{eqnarray}\label{ineqparapq}
\frac{1}{p'} \vee\frac{1}{[1/H] +1} &<& \frac{1}{p} < H,\qquad
\frac34 <\frac{1}{q} < H+ \frac12,
\nonumber
\\[-8pt]
\\[-8pt]
\nonumber
 \frac{1}{p}+\frac{1}{q}& >&1,\qquad
\frac{1}{q}-\frac{1}{p} >\frac12.
\end{eqnarray}
For example, $1/p= H -2\varepsilon$ and $1/q = H+1/2 -\varepsilon$
for sufficiently
small $\varepsilon>0$
satisfy (\ref{ineqparapq}).
Indeed,
\[
\frac{1}{p}+\frac{1}{q} = 1 +2\biggl(H -\frac14\biggr) -3\varepsilon,\qquad
\frac
{1}{q}-\frac{1}{p} =\frac12 +\varepsilon.
\]

For this $p$ and $q=\delta^{-1}$, the fBM with the Hurst parameter $H$
can be lifted to $G\Omega_p ({\mathbf{R}}^d)$
and its Cameron--Martin space $\mathcal{H}^H$ satisfies Proposition~\ref
{prFV} above.
In particular, the Young integral of $k \in\mathcal{H}^H$ with respect to
itself is possible since $q<2$.
The shift and the pairing of $X \in G\Omega_p ({\mathbf{R}}^d)$ by $k \in
\mathcal{H}^H$ can be defined
since $1/p+1/q>1$.
In what follows we always assume (\ref{ineqparapq}).
%

%%%%%%%%%%%%%%%%%%%%%
%
The Banach space $W^{\delta,p}$ is defined by (\ref{defbesovnorm}).
In Adams~\cite{ad}, its original definition is given by a kind of real
interpolation (precisely, the trace space of J. L. Lions, see paragraph
7.35,~\cite{ad}) of
$W^{1,p}$ and $W^{0,p}= L^p$.
Those are equivalent Banach spaces (paragraph 7.48,~\cite{ad}).
On the other hand, $L^{\delta,p}$ is defined by the complex
interpolation of (the complexification of)
$W^{1,p}$ and $W^{0,p}= L^p$, that is, $L^{\delta, p} = [W^{1,p},
L^p]_{1-\delta}$.
If $p=2$, $L^{\delta, 2} $ and (the complexification of) $W^{\delta
,2}$ are equivalent Hilbert spaces
(not unitarily equivalent, see paragraph 7.59,~\cite{ad}).
As a result , $L_{\mathrm{real}}^{\delta, 2} $ and $W^{\delta,2}$ are equivalent
real Hilbert spaces,
where $L_{\mathrm{real}}^{\delta, 2} $ is the subspace of ${\mathbf{R}}^d$-valued
functions in $L^{\delta, 2} $.
%

%th4.3 #&#
\begin{tm}\label{tmONB}
The following functions of $t \in[0,1]$
form an orthonormal basis of $L_{\mathrm{real}}^{\delta, 2} $ and of $L^{\delta
, 2} = L_{\mathrm{real}}^{\delta, 2} \otimes\mathbf{C} $:
\[
\{ 1\cdot\mathbf{e}_i  |  1 \le i \le d \} \cup\biggl\{ \frac{\sqrt{2} }{
(1+n^2)^{\delta/2} } \cos(n\pi t) \mathbf{e}_i  \Big| n \ge1, 1 \le i \le d
\biggr\}.
\]
Here, $\{\mathbf{e}_1, \ldots, \mathbf{e}_d\}$ is the canonical orthonormal
basis of ${\mathbf{R}}^d$.
\end{tm}

\begin{pf}
It is sufficient to prove the case $d=1$.
Note that
\[
W^{1,2} = \Biggl\{ f= c_0 + \sum_{n=1}^{\infty} c_n \sqrt{2} \cos(n\pi t)  \Big|
c_n \in\mathbf{C},
\sum_{n=0}^{\infty} (1+n^2) |c_n|^2 <\infty\Biggr\}
\]
and $\|f\|^2_{ W^{1,2}} =\sum_{n=0}^{\infty} (1+n^2) |c_n|^2$.
Similarly,
\[
L^{2} = \Biggl\{ f= c_0 + \sum_{n=1}^{\infty} c_n \sqrt{2} \cos(n\pi t)  \Big|
c_n \in\mathbf{C},
\sum_{n=0}^{\infty} |c_n|^2 <\infty\Biggr\}
\]
and $\|f\|^2_{ L^{2}} =\sum_{n=0}^{\infty} |c_n|^2$.
Therefore, $W^{1,2}$ and $L^2$ are unitarily isometric to $l_2^{(1)}$
and $l_2^{(0)}=l_2$, respectively,
where
\[
l_2^{(\delta)} =
\Biggl\{
\mathbf{c} =(c_n)_{n=0,1,2,\ldots} \in\mathbf{C}^{\infty}  \Big|
\| \mathbf{c} \|^2_{ l_2^{(\delta)} } =\sum_{n=0}^{\infty}
(1+n^2)^{\delta
} |c_n|^2\Biggr\}\qquad
 (\delta\in{\mathbf{R}}).
\]
Thus, the problem is reduced to the complex interpolation of two
Hilbert spaces of sequences.
A simple calculation shows that $[ l_2^{(1)} ,l_2 ]_{1 -\delta} =
l_2^{(\delta)} $.
This implies
\[
L^{\delta, 2} = \Biggl\{ f= c_0 + \sum_{n=1}^{\infty} c_n \sqrt{2} \cos
(n\pi
x)  \Big|  c_n \in\mathbf{C},
\sum_{n=0}^{\infty} (1+n^2)^{\delta} |c_n|^2 <\infty\Biggr\}
\]
with
$\|f\|^2_{ L^{\delta,2}} =\sum_{n=0}^{\infty} (1+n^2)^{\delta} |c_n|^2$,
which ends the proof.
\end{pf}

%%%%%%%%%%%%%%%%%%%%%%%%%%%%%%%%%%%%%

%
We compute the $p$-variation norm of cosine functions.
The following lemma is taken from Nate Eldredge's unpublished
manuscripts~\cite{eld}.
Before stating it, we introduce some definitions.
Let $x$ be a one-dimensional continuous path with $x_0=0$. We say that
$s \in[0,1]$ is a forward maximum
(or forward minimum)
if $x_s =\max x|_{[s,1]}$ (or $x_s =\min x|_{[s,1]}$, resp.).
Suppose $x$ is piecewise monotone with local extrema $\{ 0=s_0 <s_1<s_2
<\cdots< s_n=1 \}$.
(For simplicity, we assume $s_0,s_2, \ldots$ are local minima and $s_1,
s_3, \ldots$ are local maxima.
The reverse case is easily dealt with by just replacing $x$ with $-x$.)
If $s_2, s_4, \ldots$ are not only local minima but also forward minima,
and $s_1, s_3, \ldots$ are not only local maxima but also forward maxima,
then we say $x$ is jog-free.
(Note that $x_0$ is not required to be a forward extremum.)
%
%pr4.4 #&#
\begin{pr}\label{preld}
Let $p \ge1$.
\textup{(i)}
If a one-dimensional continuous path $x$ with $x_0=0$ is jog-free with
extrema $\{ 0=s_0 <s_1<s_2 <\cdots< s_n=1 \}$,
then
\[
\|x\|_{p\mbox{-}\mathrm{var}} = \Biggl(
\sum_{i=1}^n | x_{s_i} -x_{s_{i-1}}|^p
\Biggr)^{1/p}.
\]\vspace*{-6pt}

\begin{longlist}
\item[(ii)]
In particular, the $p$-variation norm of $c_n(t)= \cos( n \pi t) -1$ is
given by
$
\|c_n \|_{p\mbox{-}\mathrm{var}} = 2 n^{1/p}.
$\vadjust{\goodbreak}
\end{longlist}
\end{pr}

\begin{pf}
(ii) is immediate from \textup{(i)}. We show \textup{(i)}.
For a continuous path $y$ and a partition $\mathcal{P}=\{ 0=t_0 <t_1<t_2
<\cdots< t_n=1 \}$,
we set
$
V_{p, \mathcal{P}} (y) = ( \sum_{i=1}^n
| y_{t_i} - y_{t_{i-1}} |^p
)^{1/p}
$.
Then, $ \|y \|_{p\mbox{-}\mathrm{var}} = \sup_{\mathcal{P}} V_{p, \mathcal{P}} (y)$.
First, note that if $y$ is monotone increasing (or decreasing) on $[
t_{i-1}, t_{i+1} ]$, then it is easy to see that
$ V_{p, \mathcal{P} \setminus\{ t_i\}} (y) \ge V_{p, \mathcal{P}} (y)$.
In other words, intermediate points in monotone intervals should not be
included.

Let $x$ be jog-free with extrema $\mathcal{Q}= \{ 0=s_0 <s_1<s_2 <\cdots<
s_n=1 \}$
as in the statement of \textup{(i)} and let $\mathcal{P}=\{ 0=t_0 <t_1<t_2
<\cdots< t_n=1 \}$
be a partition which does not include all the $s_j$'s.
We will show below that there exists an $s_j$ such that $V_{p, \mathcal{P}
\cup\{ s_j\}} (x) \ge V_{p, \mathcal{P}} (x)$.

Let $s_j$ be the first extremum not contained in $\mathcal{P}$.
(For simplicity, we assume it is local and forward maximum.)
Let $t_i$ be the last element of $\mathcal{P}$ less than $s_j$. Then,
$s_{j-1} \le t_i \le s_j \le t_{i+1}$.
Since $x$ is increasing on $[s_{j-1} , s_j]$ and $x_{s_j}$ is forward maximum,
\[
x_{s_j}-x_{t_i} \ge x_{t_{i+1}} -x_{t_i},\qquad
x_{s_j}-x_{t_{i+1}} \ge x_{t_{i}} -x_{t_{i+1}},
\]
which yields that
$
| x_{s_j}-x_{t_i} |^p +| x_{s_j}-x_{t_{i+1}}|^p \ge|x_{t_{i+1}}
-x_{t_i}|^p.$
Therefore,\break $V_{p, \mathcal{P} \cup\{ s_j\}} (x) \ge V_{p, \mathcal{P}} (x)$.

For any $\varepsilon>0$, there exists $\mathcal{P}$ such that $V_{p,
\mathcal{P}}
(x) \ge\|x \|_{p\mbox{-}\mathrm{var}} -\varepsilon$.
First by adding all the $s_j$'s, then by removing all the intermediate
points (i.e., $t_i$'s which are not one of $s_j$'s), we get
$V_{p, \mathcal{Q}} (x) \ge\|x \|_{p\mbox{-}\mathrm{var}} -\varepsilon$.
Letting $\varepsilon\searrow0$, we complete the proof of \textup{(i)}.
\end{pf}

%%%%%%%%%%%%%%%%%%%%%%%%%%%%%%%
%
Now we calculate the Hessian of $\Psi$, which is defined in (\ref{rdeqvar}).
For $q<2$, ODE like (\ref{rdeqvar}) is well-defined in the
$q$-variation sense, thanks to the Young integral.
The continuity of $\Psi$ is well known. Smoothness of the It\^o map in
the $q(<2)$-variation setting is studied in Li--Lyons~\cite{ll}.
The explicit form of the derivatives are obtained in a similar way to
the case of (stochastic) Taylor expansion.

Let $q \in[1,2)$ for a while and fix $\gamma\in C^{q\mbox{-}\mathrm{var}}_0$.
Then $\phi^0 = \Psi(\gamma)$ is also of finite $q$-variation, which
takes values in ${\mathbf{R}}^n$.
Set
\[
d \Omega_t =
\nabla\sigma(\phi^0_t )\langle \cdot  , d\gamma_t \rangle+
\nabla\beta_0 (\phi^0_t) \langle \cdot  \rangle \,dt.
\]
Then, $\Omega$ is an $\operatorname{End }({\mathbf{R}}^n)$-valued path of finite
$q$-variation.
Next, consider the following $\operatorname{End }({\mathbf{R}}^n)$-valued ODE in the
$q$-variation sense:
%
%e4.3 #&#
\begin{equation}\label{eqeqM}
dM_t = d \Omega_t \cdot M_t, \qquad  M_0=\mathrm{Id}_n.
\end{equation}
Its inverse satisfies a similar ODE:
%
%e4.4 #&#
\begin{equation}\label{eqeqM-1}
dM_t^{-1} = -M_t^{-1} \cdot d \Omega_t,\qquad   M_0^{-1}=\mathrm{Id}_n,
\end{equation}
although the coefficients of these ODEs are not bounded, thanks to
their special forms, to
a unique solution
[For this kind of equation with unbounded coefficients, existence of a
local solution and uniqueness are easier.
The problem is existence of a global solution.
If $M$ is a local solution of (\ref{eqeqM}) and $ a \in\operatorname{GL}(n
,{\mathbf{R}})$, then\vadjust{\goodbreak} $Ma$ is a local solution to (\ref{eqeqM})
with a initial condition $M_0=a$.
This fact, combined with existence of a local solution, implies
existence of a global solution.]
The map $\gamma\mapsto M$ in the $q$-variational setting is locally
Lipschitz continuous.
(In this paper, however, $\gamma$ is always fixed and, hence, so are
$\Omega$ and $M$.)
If $\gamma$ is controlled by a control function $\omega$,
then $M$ and $M^{-1}$ are controlled by
$\hat\omega(s,t) =C(\omega(s,t) +(t-s))$, where $C>0$ is a constant
which depends on $q$ and $\omega(0,1)$.
A rigorous proof for this paragraph can be found in~\cite{ina2}, for instance.

%%%%%%%%%%%%%%%%%%%%%%%%%%%%%%%%%%%%%%%%%%

Set $\chi(k) = (\nabla\Psi)(\gamma) \langle k \rangle$ for simplicity.
This is
a continuous path of
finite $q$-variation, if $k$ is of finite $q$-variation.
Then, it satisfies an ${\mathbf{R}}^n$-valued ODE:
%
%e4.5 #&#
\begin{equation}\label{defchik}
d \chi_t - \nabla\sigma(\phi^0_t) \langle\chi_t, d\gamma_t
\rangle-
\nabla
\beta_0 (\phi^0_t) \langle\chi_t \rangle \,dt
=
\sigma(\phi^0_t) \,dk_t,\qquad
 \chi_0=0.
\end{equation}
From this, we can obtain an explicit expression as follows:
%
%e4.6 #&#
\begin{equation}\label{defchik2}
\chi(k)_t = (\nabla\Psi)(\gamma) \langle k \rangle_t = M_t \int
_0^t M^{-1}_s
\sigma( \phi^0_s) \,dk_s.
\end{equation}
Note that the right-hand side is a Young integral and $k \mapsto\chi
(k)$ extends to a continuous map
from $C^{p \mbox{-}\mathrm{var}}_0 ({\mathbf{R}}^d)$ to $C^{p \mbox{-}\mathrm{var}}_0 ({\mathbf{R}}^n)$.

%%%%%%

In a similar way, $\psi_t =\nabla^2 \Psi(\gamma) \langle k,k
\rangle_t$
satisfies the following ODE:
%
%e4.7 #&#
\begin{eqnarray}\label{defpsi}
&&
d\psi_t - \nabla\sigma(\phi^0_t) \langle\psi_t, d\gamma_t
\rangle-
\nabla
\beta_0 (\phi^0_t) \langle\psi_t \rangle \,dt
\nonumber
\\
&&\qquad=
2 \nabla\sigma(\phi^0_t) \langle\chi(k)_t, dk_t \rangle
+ \nabla^2 \sigma(\phi^0_t) \langle\chi(k)_t, \chi(k)_t ,d\gamma_t
\rangle
\\
&&\qquad\quad{}+ \nabla^2 \beta_0 (\phi^0_t) \langle\chi(k)_t, \chi(k)_t \rangle \,dt,
\qquad
\psi_0=0.\nonumber
\end{eqnarray}
From this and by polarization, we see that
%
%e4.8 #&#
\begin{eqnarray}\label{eqhessexp}
&&
\nabla^2 \Psi(\gamma) \langle f, k \rangle_t\nonumber\\[-2pt]
&&\qquad=
M_t \int_0^t M^{-1}_s
\{
\nabla\sigma(\phi^0_s) \langle\chi(f)_s, dk_s \rangle
+
\nabla\sigma(\phi^0_s) \langle\chi(k)_s, df_s \rangle
\}
\nonumber\\[-2pt]
&&\qquad\quad{}
+
M_t \int_0^t M^{-1}_s
\{
\nabla^2 \sigma(\phi^0_t) \langle\chi(f)_s, \chi(k)_s ,d\gamma_s
\rangle
\\[-2pt]
&&\hspace*{80pt}\qquad{}+\nabla^2 \beta_0 (\phi^0_s) \langle\chi(f)_s, \chi(k)_s \rangle \,ds
\}
\nonumber\\[-2pt]
&&\qquad=:  V_1 (f,k)_t + V_2 (f,k)_t.\nonumber
\end{eqnarray}
It is obvious that
\[
(f,k) \in C^{q\mbox{-}\mathrm{var}}_0 ({\mathbf{R}}^d) \times C^{q\mbox{-}\mathrm{var}}_0 ({\mathbf{R}}^d)
\mapsto\nabla^2 \Psi(\gamma) \langle f, k \rangle\in C_0^{q\mbox{-}\mathrm{var}}
({\mathbf{R}}^n)
\]
is a symmetric bounded bilinear functional.

%%%%%%%%

Note that $\chi(k)$ and $\psi(k,k)/2$ are similar to $\phi^1$ and
$\phi^2$, respectively, when $X =k$.
Indeed, they are the first and the second term in the Taylor expansion
for $\hat\Phi_0 ( \varepsilon X +\gamma, \lambda)$.
[See (\ref{defphi1}) and (\ref{eqodephi2}) below and compare.]\vspace*{1pt}
Therefore, $k \mapsto\chi(k), \psi(k,k)$ extend to continuous maps
from $G\Omega_p ({\mathbf{R}}^d)$ to $ C^{p\mbox{-}\mathrm{var}}_0 ({\mathbf{R}}^n)$.
We will write $ \chi(X), \psi(X,X)$ for $X \in G\Omega_p ({\mathbf{R}}^d)$.
%%%%%%%%%%%%

%le4.5 #&#
\begin{lm} \label{lmtraceV2}
Let $1/4 <H< 1/2$ and choose $p$ and $q$ as in (\ref{ineqparapq})
Then, for any bounded linear functional\vadjust{\goodbreak} $\alpha\in C^{p\mbox{-}\mathrm{var}}_0({\mathbf{R}}^n)^*$,\vspace*{-1pt}
the symmetric bounded bilinear form $\alpha\circ V_2 \langle \cdot
,
\cdot \rangle$ on
the Cameron--Martin space $\mathcal{H}^H$ is of trace class.
In particular, if $p' \ge p$, $\nabla F (\phi^0) \circ V_2$ is of
trace class
for a Fr\'echet differentiable function $F\dvtx  C^{p'\mbox{-}\mathrm{var}}_0({\mathbf{R}}^n) \to
{\mathbf{R}}$.
Moreover, $\alpha\circ V_2$ extends to a bounded bilinear form on
$C^{p\mbox{-}\mathrm{var}}_0 ({\mathbf{R}}^d)$.
A similar fact holds for $\nabla^2 F (\phi^0) \langle\chi( \cdot ),
\chi
( \cdot )\rangle$, too.\looseness=-1\vspace*{-3pt}
\end{lm}

\begin{pf}
Since $t \mapsto M_t$ and $t \mapsto M_t^{-1} \sigma(\phi^0_t)$ are of
finite $q$-variation,
the map
$h \mapsto\chi(h)$ extends to a bounded linear map from $C_0^{p\mbox{-}\mathrm{var}}
({\mathbf{R}}^d)$
to $C_0^{p\mbox{-}\mathrm{var}} ({\mathbf{R}}^n)$, thanks to the Young integral.
By using the Young\vspace*{1pt} integral again,
we see that
$(h,k) \mapsto V_2(h,k)$ extends to a bounded bilinear map from
$C_0^{p\mbox{-}\mathrm{var}} ({\mathbf{R}}^d) \times C_0^{p\mbox{-}\mathrm{var}} ({\mathbf{R}}^d)$
to $C_0^{q\mbox{-}\mathrm{var}} ({\mathbf{R}}^n) \subset C_0 ({\mathbf{R}}^d)$.

On the other hand, $\mu^H$ (the law of the fBM with the Hurst parameter
$H$) is supported in $C_0^{p\mbox{-}\mathrm{var}} ({\mathbf{R}}^d)$.
In other words,
$(
\mathcal{X}, \mathcal{H}^H, \mu^H)$ is an abstract Wiener space,
where
$\mathcal{X}$ is the closure of $\mathcal{H}^H$ with respect to the
$p$-variation norm.
[According to Jain and Monrad~\cite{jain}, pages 47--48,
$C_0^{p\mbox{-}\mathrm{var}} ({\mathbf{R}}^d)$ is not separable and, consequently, $\mathcal{H}^H$ cannot be dense in $C_0^{p\mbox{-}\mathrm{var}} ({\mathbf{R}}^d)$.
So, we use
$\mathcal{X}$ instead of $C_0^{p\mbox{-}\mathrm{var}} ({\mathbf{R}}^d)$, because an abstract
Wiener space must be separable by definition.]

Therefore, $\alpha\circ V_2$ is a bounded bilinear form on an abstract
Wiener space.
By Goodman's theorem (Theorem 4.6, Kuo~\cite{kuo}), its restriction on
the Cameron--Martin space is of trace class.\vspace*{-3pt}
\end{pf}

Now we compute $V_1$.\vspace*{-3pt}
%
%le4.6 #&#
\begin{lm} \label{lmhilschV1}
Let $1/4 <H< 1/2$ and choose $p$ and $q$ as in (\ref{ineqparapq}).
Then, for any bounded linear functional $\alpha\in C^{p\mbox{-}\mathrm{var}}_0({\mathbf{R}}^n)^*$,
the symmetric bounded bilinear form $\alpha\circ V_1 \langle \cdot
,
\cdot \rangle$ on
the Cameron--Martin space $\mathcal{H}^H$ is Hilbert--Schmidt.
In particular, if $p' \ge p$, $\nabla F (\phi^0) \circ V_1$ is Hilbert--Schmidt
for a Fr\'echet differentiable function $F\dvtx  C^{p'\mbox{-}\mathrm{var}}_0({\mathbf{R}}^n) \to
{\mathbf{R}}$.
Moreover, if $\alpha_l$ is weak* convergent to $\alpha$ as $l \to
\infty
$ in $C^{p\mbox{-}\mathrm{var}}_0({\mathbf{R}}^n)^*$,
then $\alpha_l \circ V_1$ converges to $\alpha\circ V_1$ as $l \to
\infty$ in the Hilbert--Schmidt norm.\vspace*{-3pt}
\end{lm}

The rest of this section is devoted to proving this lemma.
An integration by parts yields that
\[
V_1 \langle f,k\rangle=R_1 \langle f,k\rangle+ R_1 \langle k,
f\rangle- ( R_2 \langle f,k\rangle
+R_2 \langle k,f\rangle),
\]
where, from (\ref{defchik2}),
\begin{eqnarray*}
R_1 \langle f,k\rangle_t
&=&
M_t \int_0^t M^{-1}_s \nabla\sigma(\phi^0_s)
\langle\sigma(\phi^0_s) f_s ,dk_s \rangle,
\\[-2pt]
R_2 \langle f,k\rangle_t
&=&
M_t \int_0^t M^{-1}_s \nabla\sigma(\phi^0_s)
\biggl\langle M_s \int_0^s d[M_u^{-1} \sigma(\phi^0_u) ] f_u ,dk_s
\biggr\rangle.\vspace*{-3pt}
\end{eqnarray*}

%le4.7 #&#
\begin{lm}\label{lmhsK2}
Let $R_2$ be as above and $\alpha\in C^{p\mbox{-}\mathrm{var}}_0({\mathbf{R}}^n)^*$.
Then, as a bilinear form on $\mathcal{H}^H$,
$\alpha\circ R_2$ is of trace class.
Moreover,\vadjust{\goodbreak} if $\alpha_l$ is weak* convergent to $\alpha$ as $l \to
\infty
$ in $C^{p\mbox{-}\mathrm{var}}_0({\mathbf{R}}^n)^*$,
then $\alpha_l \circ R_2$ converges to $\alpha\circ R_2$ as $l \to
\infty$ in the Hilbert--Schmidt norm.
\end{lm}

\begin{pf}
We use the Young integral.
Since $u \mapsto M_u^{-1} \sigma(\phi^0_u) $ is of finite $q$-variation,
we see that
\[
\biggl\| \int_0^{\cdot} d[M_u^{-1} \sigma(\phi^0_u) ] f_u \biggr\|_{q\mbox{-}\mathrm{var}}
\le
c_1 \| M_{\cdot}^{-1} \sigma(\phi^0_{\cdot}) \|_{q\mbox{-}\mathrm{var}} \| f \|_{p\mbox{-}\mathrm{var}}
\le
c_2 \| f \|_{p\mbox{-}\mathrm{var}}.
\]
Similarly, since
$s \mapsto M_s^{-1} \nabla\sigma(\phi^0_s), M_s$ are of finite $q$-variation,
%
%e4.9 #&#
\begin{eqnarray}\label{ineqhsK2}
&&\biggl\|
M_{\cdot}
\int_0^{\cdot} M^{-1}_s \nabla\sigma(\phi^0_s)
\biggl\langle M_s \int_0^s d[M_u^{-1} \sigma(\phi^0_u) ] f_u ,dk_s
\biggr\rangle
\biggr\|_{p\mbox{-}\mathrm{var}}
\nonumber
\\[-8pt]
\\[-8pt]
\nonumber
&&\qquad
\le
c_3 \| f \|_{p\mbox{-}\mathrm{var}} \|k\|_{p\mbox{-}\mathrm{var}}.
\end{eqnarray}
Thus, $(f ,k) \mapsto R_2 \langle f,k\rangle$ is a bounded bilinear
map from
$C_0^{p\mbox{-}\mathrm{var}}({\mathbf{R}}^d) \times C_0^{p\mbox{-}\mathrm{var}}({\mathbf{R}}^d) $ to
$C_0^{p\mbox{-}\mathrm{var}}({\mathbf{R}}^n) \subset C_0 ({\mathbf{R}}^n)$.
In particular,
$\alpha\circ R_2$ is a bounded bilinear\vspace*{1pt} form on $C_0^{p\mbox{-}\mathrm{var}}({\mathbf{R}}^n)$.
Again, by Goodman's theorem (Theorem 4.6,~\cite{kuo}), its restriction
on the Cameron--Martin space is of trace class.

Now we prove the convergence. Note that (\ref{ineqhsK2}) still holds
even when $f$ or $k$
do not start at $0$.
Consider the following continuous inclusions
(see Proposition~\ref{prFV}. Below, all the function space is ${\mathbf{R}}^d$-valued):
\[
\mathcal{H}^H \hookrightarrow W_0^{\delta, 2} \cong L^{\delta,
2}_{0,\mathrm{real}} \hookrightarrow L^{\delta, 2} _{\mathrm{real}}
\hookrightarrow C^{q\mbox{-}\mathrm{var}} \hookrightarrow C^{p\mbox{-}\mathrm{var}},
\]
where $\delta=1/q$ and $\cong$ denotes isomorphism (but not unitary)
of Hilbert spaces.
Let us first consider $R_2 \vert_{L^{\delta, 2}_{\mathrm{real}} \times
L^{\delta
, 2}_{\mathrm{real}} }$.
We will show that,\vspace*{-1pt} for an ONB $\{f_k\}_{k=1,2,\ldots}$ of $L^{\delta,
2}_{\mathrm{real}}$, it holds that
$
\sum_{k,j=1}^{\infty} \| R_2 \langle f_k,f_j \rangle\|^2_{p\mbox{-}\mathrm{var}}
<\infty.
$
As in Theorem~\ref{tmONB}, we set
$f_{0,i} (t) =1\cdot\mathbf{e}_i$ and $f_{m,i} (t)= (1+m^2)^{-\delta/2}
\sqrt{2} \cos(m\pi t) \mathbf{e}_i \ (m=1,2,\ldots).$
By Proposition~\ref{preld},
\[
\| f_{m,i} \|_{p\mbox{-}\mathrm{var}} \le(1+m^2)^{-\delta/2} \sqrt{2} (1+ 2m^{1/p})
\le
c \biggl( \frac{1}{ 1+m} \biggr)^{{1}/{q}- {1}/{p}}
\]
for some constant $c>0$.
From this and (\ref{ineqhsK2}),
\begin{eqnarray*}
&&\sum_{i,i' =1}^d \sum_{m,m' =0}^{\infty} \| R_2 \langle f_{m,i},f_{m',i'}
\rangle\|^2_{p\mbox{-}\mathrm{var}}\\
&&\qquad\le
c \sum_{i,i' =1}^d \sum_{m,m' =0}^{\infty} \| f_{m,i} \|_{p\mbox{-}\mathrm{var}}^2
\| f_{m',i'} \|_{p\mbox{-}\mathrm{var}}^2
\\
&&\qquad\le
c \sum_{m=0}^{\infty}
\biggl( \frac{1}{ 1+m} \biggr)^{2 ({1}/{q}- {1}/{p})}
\sum_{m' =0}^{\infty}
\biggl( \frac{1}{ 1+m'} \biggr)^{2 ({1}/{q}- {1}/{p})}
<\infty,
\end{eqnarray*}
because $1/q-1/p >1/2$.
Here, the constant $c>0$ may change from line to line.

By the Banach--Steinhaus theorem,
$\| \alpha_l -\alpha\|_{C^{p\mbox{-}\mathrm{var},*}} \le c$ for some constant $c>0$.
Hence,
\[
| (\alpha_l-\alpha) \circ R_2 \langle f_{m,i},f_{m',i'} \rangle|^2
\le
c^2
\| R_2 \langle f_{m,i},f_{m',i'} \rangle\|^2_{p\mbox{-}\mathrm{var}}.
\]
By the dominated convergence theorem, $\| \alpha_k \circ R_2-\alpha
\circ R_2 \|_{\operatorname{HS}-L^{\delta,2}_{\mathrm{real}}} \to0$
as $k \to\infty$.
(The norm denotes the Hilbert--Schmidt norm.)
This implies that
\[
\| \alpha_l \circ R_2-\alpha\circ R_2 \|_{\operatorname{HS}-\mathcal{H}^H}
\le
\| \iota\|_{\operatorname{op}} \| \iota^* \|_{\operatorname{op}}
\| \alpha_l \circ R_2-\alpha\circ R_2 \|_{\operatorname{HS}-L^{\delta,2}_{\mathrm{real}}}
\to0
\]
as $l \to\infty$,
where $\iota\dvtx  \mathcal{H}^H \hookrightarrow L^{\delta,2}_{\mathrm{real}}$ denotes
the inclusion.
\end{pf}

%%%%%%%%%%%%%%%%%%%%%%%%%
%
%le4.8 #&#
\begin{lm}\label{lmhsK1}
Let $R_1$ be as above and $\alpha\in C_0^{p\mbox{-}\mathrm{var}}({\mathbf{R}}^n)^*$.
Then, as a bilinear form on $\mathcal{H}^H$,
$\alpha\circ R_1$ is Hilbert--Schmidt.
Moreover, if $\alpha_l$ is weak* convergent to~$\alpha$ as $l \to
\infty
$ in $C_0^{p\mbox{-}\mathrm{var}}({\mathbf{R}}^n)^*$,
then $\alpha_l \circ R_1$ converges to $\alpha\circ R_1$ as $l \to
\infty$ in the Hilbert--Schmidt norm.
\end{lm}

\begin{pf}
The proof is similar to the one for Lemma~\ref{lmhsK2}.
It is sufficient to show that
%
%e4.10 #&#
\begin{equation}\label{estjuu}
\sum_{i,i' =1}^d \sum_{m,m' =0}^{\infty} \| R_1 \langle f_{m,i},f_{m',i'}
\rangle\|^2_{p\mbox{-}\mathrm{var}}
<\infty.
\end{equation}
In this proof, $c>0$ is a constant which may change from line to line.

It is easy to see that, if $m \neq m'$,
\begin{eqnarray*}
\sqrt{2} \cos(m \pi t) d\bigl[ \sqrt{2} \cos(m' \pi t) \bigr]
&=&
-2m' \pi\cos(m \pi t) \sin(m' \pi t) \,dt
\\
&=&
-m' \pi\bigl\{ \sin\bigl((m'+m) \pi t\bigr) + \sin\bigl( (m'-m) \pi t\bigr) \bigr\}\,dt
\\
&=&
m' d \biggl[
\frac{ \cos((m'+m) \pi t) }{m' +m} + \frac{ \cos((m' -m) \pi t) }{m'
-m} \biggr],
\end{eqnarray*}
and that, if $m =m'$,
$
\sqrt{2} \cos(m \pi t) d[ \sqrt{2} \cos(m\pi t) ]
= d [ \cos(2m\pi t)] /2.
$

In the following, fix $i,i'$. First, we consider the case $m=m'$:
\begin{eqnarray*}
&&R_1 \langle f_{m, i}, f_{m,i'} \rangle_t\\
&&\qquad=
M_t \int_0^t M^{-1}_s \nabla\sigma(\phi^0_s)
\langle\sigma(\phi^0_s) \mathbf{e}_i , \mathbf{e}_{i'} \rangle\sqrt{2}
\frac{ \cos(m \pi s )} { (1+m^2)^{1/2q} } \, d \biggl[ \frac{ \sqrt{2}
\cos(m\pi s )}{ (1+m^2)^{1/2q}} \biggr]
\\
&&\qquad=
\frac{1/2}{(1+m^2)^{1/q}}
M_t \int_0^t M^{-1}_s \nabla\sigma(\phi^0_s)
\langle\sigma(\phi^0_s) \mathbf{e}_i , \mathbf{e}_{i'} \rangle\  d [ \cos
(2m\pi
s ) ].
\end{eqnarray*}
By the Young integral and Proposition~\ref{preld}, we see that
\begin{eqnarray*}
\| R_1 \langle f_{m, i}, f_{m,i'} \rangle\|_{p\mbox{-}\mathrm{var}}^2
&\le&
\frac{c}{(1+m^2)^{2/q}} \| \cos( 2m \pi \cdot ) -1 \|_{p\mbox{-}\mathrm{var}}^2
\\
&\le&
\frac{c m^{2/p} }{(1+m^2)^{2/q}}
\le
\frac{c}{(1+m)^{4/q - 2/p} }.
\end{eqnarray*}
Since $4/q - 2/p >1$,
%
%e4.11 #&#
\begin{equation} \label{estsono1}
\sum_{m=0}^{\infty}
\| R_1 \langle f_{m, i}, f_{m,i'} \rangle\|_{p\mbox{-}\mathrm{var}}^2 <\infty.
\end{equation}

Next we consider the case $m\neq m'$:
\begin{eqnarray*}
&&
R_1 \langle f_{m, i}, f_{m' , i'} \rangle_t
\\
&&\qquad=
M_t \int_0^t M^{-1}_s \nabla\sigma(\phi^0_s)
\langle\sigma(\phi^0_s) \mathbf{e}_i , \mathbf{e}_{i'} \rangle\sqrt{2}
\frac{ \cos(m \pi s)} { (1+m^2)^{1/2q} } \, d \biggl[ \frac{ \sqrt{2}
\cos
(m' \pi s)}{ (1+m'^2)^{1/2q}} \biggr]\\
&&\qquad=
\frac{m'}{ (1+m^2)^{1/2q} (1+m'^2)^{1/2q} (m' +m)}\\
&&\qquad\quad{}
  \times
M_t \int_0^t M^{-1}_s \nabla\sigma(\phi^0_s)
\langle\sigma(\phi^0_s) \mathbf{e}_i , \mathbf{e}_{i'} \rangle
\, d \bigl[ \cos\bigl((m'+m) \pi s \bigr) \bigr]\\
&&\qquad\quad{}+
\frac{m'}{ (1+m^2)^{1/2q} (1+m'^2)^{1/2q} (m' -m)}\\
&&\qquad\qquad{}
  \times
M_t \int_0^t M^{-1}_s \nabla\sigma(\phi^0_s)
\langle\sigma(\phi^0_s) \mathbf{e}_i , \mathbf{e}_{i'} \rangle
\, d \bigl[ \cos\bigl((m'-m) \pi s \bigr) \bigr]\\
&&\qquad=:
\hat{R}_1^{i,i'} (m,m')_t +\hat{R}_2^{i,i'} (m,m')_t.
\end{eqnarray*}

By using the estimate for the Young integral again, we see that
\begin{eqnarray*}
&&\| \hat{R}_1^{i,i'} \langle f_{m, i}, f_{m,i'}
\rangle\|_{p\mbox{-}\mathrm{var}}^2\\
&&\qquad\le
\frac{c m'^2 }{(1+m^2)^{1/q} (1+m'^2)^{1/q} |m'+m|^2} \bigl\| \cos\bigl( (m' +m)
\pi \cdot \bigr) -1 \bigr\|_{p\mbox{-}\mathrm{var}}^2
\\
&&\qquad\le
\frac{c m'^2 |m'+m|^{2/p} }{(1+m^2)^{1/q} (1+m'^2)^{1/q} |m'+m|^2}
\\
&&\qquad\le
\frac{c |(m' +m) -m |^{ 2 (1 -1/q)} }{(1+|m|)^{2/q} |m'+m|^{2 (1-1/p)
} }
\\
&&\qquad\le
\frac{c}{ (1+|m|)^{2/q} (1+ |m'+m|)^{2 (1/q-1/p) } }
\\
&&\qquad\quad
{}
+
\frac{c}{ (1+|m|)^{4(1/q -1/2)} (1+ |m'+m|)^{2 (1-1/p) } }.
\end{eqnarray*}
It is easy to see that $2/q >1$ and $2(1-1/p) >1$ hold.
From (\ref{ineqparapq}), $2 (1/q-1/p) >1$ and $4(1/q -1/2)>1$.
(The condition $1/q >3/4$ is used here.)
Therefore,
%
%e4.12 #&#
\begin{eqnarray}\label{estsono2}
&&
\sum_{ 0 \le m,m' < \infty, m\neq m'}
\| \hat{R}_1^{i,i'} \langle f_{m, i}, f_{m,i'} \rangle\|_{p\mbox{-}\mathrm{var}}^2
\nonumber\\
&&\qquad\le
c \sum_{m, m' \in\mathbf{Z} }
\biggl(
\frac{1}{ (1+|m|)^{2/q} (1+ |m'+m|)^{2 (1/q-1/p) } }
\nonumber\\
&&\hspace*{41pt}\qquad\quad{}
+
\frac{1}{ (1+|m|)^{4(1/q -1/2)} (1+ |m'+m|)^{2 (1-1/p) } }\biggr)
\\
&&\qquad=
c \sum_{m \in\mathbf{Z} }
\frac{1}{ (1+|m|)^{2/q} }
\biggl(
\sum_{m' \in\mathbf{Z}}
\frac{1}{(1+ |m'+m|)^{2 (1/q-1/p) } }
\biggr)
\nonumber\\
&&\qquad\quad{}
+
c \sum_{m \in\mathbf{Z} }
\frac{1}{ (1+|m|)^{4(1/q-1/2)} }
\biggl(
\sum_{m' \in\mathbf{Z}}
\frac{1}{(1+ |m'+m|)^{2 (1-1/p) } }
\biggr)
<\infty.\nonumber
\end{eqnarray}
In the same way as above,
%
%e4.13 #&#
\begin{equation}\label{estsono3}
\sum_{ 0 \le m,m' < \infty, m\neq m'}
\| \hat{R}_2^{i,i'} \langle f_{m, i}, f_{m,i'} \rangle\|_{p\mbox{-}\mathrm{var}}^2
<\infty.
\end{equation}
From (\ref{estsono1}), (\ref{estsono2}) and (\ref{estsono3}), we
have (\ref{estjuu}),
which completes the proof.~%
\end{pf}

%%%%%%%%%%%%%%%%%%%%%%%%%%%%%%%%%%%%%%%%%%%%%%%%%%%%%%%%%%
%%%%%%%%%%%%%%%%%%%%%%%%%%%%%%%%%%%%%%%%%%%%%%%%%%%%%%%%%%%
%%%%%%%%%%%%%%%%%%%%%%%%%%%%%%%%%%%%%%%%%%%%%%%%%%%%%%%%%%
%%%%%%%%%%%%%%%%%%%%%%%%%%%%%%%%%%%%%%%%%%%%%%%%%%%%%%%%%%%
%%%%%%%%%%%%%%%%%%%%%%%%%%%%%%%%%%%%%%%%%%%%%%%%%%%%%%%%%%
%%%%%%%%%%%%%%%%%%%%%%%%%%%%%%%%%%%%%%%%%%%%%%%%%%%%%%%%%%%
%s5 #&#
\section{A probabilistic representation of Hessian}\label{sec5}
%
%%%%%%%%%%%%%%%%%%%%%%%%%%%%%%%%%%%%%%%%%%%%%%%%%%%%%%%%%%%
%%%%%%%%%%%%%%%%%%%%%%%%%%%%%%%%%%%%%%%%%%%%%%%%%%%%%%%%%%%%
Throughout this section we assume (\ref{ineqparapq}).
Let $(\mathcal{X}, \mathcal{H}^{H}, \mu^{H})$ be the abstract Wiener space
for the fBM as in the previous section.
Here, $\mathcal{X}$ is the closure of the Cameron--Martin\vspace*{1pt} space $\mathcal{H}^{H}$ in $C_0^{p \mbox{-}\mathrm{var}} ({\mathbf{R}}^d)$.
A generic element of $\mathcal{X}$ is denoted by $w^H$.
Under $\mu^{H}$, $(w^H_t)_{0 \le t \le1}$ is the canonical realization
of $d$-dimensional fBM.

Any $\langle k,   \cdot \rangle\in( \mathcal{H}^{H})^*$ extends to a
measurable linear functional on $\mathcal{X}$,
which is denoted by $\langle k, w^H \rangle$ with a slight abuse of notation.
It satisfies
\[
\int_{\mathcal{X}} e^{ \sqrt{-1} \langle k, w^H \rangle} \mu^{H} (d
w^H) = e^{
\| k\|^2_{ \mathcal{H}^{H}} /2}.
\]
For a cylinder function $F(w^H)= f( \langle k_1, w^H \rangle, \ldots,
\langle k_m,
w^H \rangle)$, where $f \dvtx  {\mathbf{R}}^m \to{\mathbf{R}}$
is a bounded smooth function with bounded derivatives,
we set
\[
D_{k} F(w^H) = \sum_{j=1}^m \partial_j f ( \langle k_1, b \rangle,
\ldots,
\langle
k_m, b \rangle) ( k_j, k )_{\mathcal{H}^H},
\qquad
k \in\mathcal{H}^H,
\]
and
\[
D F( w^H) = \sum_{j=1}^m \partial_j f ( \langle k_1, w^H \rangle,
\ldots,
\langle
k_m, w^H\rangle) k_j.
\]
Note that $DF$ is an $\mathcal{H}^H$-valued function.\vadjust{\goodbreak}

%%%%%%%%%%%%

Let
$\mathcal{C}_n = \mathcal{C}_n(\mu^H)   (n=0,1,2,\ldots)$ be the $n$th Wiener
chaos of $w^H$.
It is well known that $\mathcal{C}_n$ are mutually orthogonal and $L^2(\mu
^H) = \bigoplus_{n=0}^{\infty} \mathcal{C}_n$.
For example, $\mathcal{C}_0 =\{ \mbox{constants} \}$
and $\mathcal{C}_1 = \{ \langle k,   \cdot \rangle |  k \in\mathcal{H}^H\}$.
The second Wiener chaos $\mathcal{C}_2$ is unitarily isometric with the
space of
symmetric Hilbert--Schmidt operators (or symmetric Hilbert--Schmidt
bilinear forms)
$\mathcal{H}^H \otimes_{sym} \mathcal{H}^H$
in a natural way.

%%%%%%%%%%%%%%%%%

%le5.1 #&#
\begin{lm}\label{lmclos}
Let $V_1$ be as in (\ref{eqhessexp}) and consider $V_1 (w^H(m), w^H (m))_t$,
where $w^H(m)$ denotes the $m$th dyadic polygonal approximation of $w^H$.
Then, for $k, \hat{k} \in\mathcal{H}^H$,
\begin{eqnarray*}
\tfrac12 D_k V_1 (w^H(m), w^H (m))_t
&=&
V_1 ( k (m), w^H (m))_t,
\\
\tfrac12 D_{\hat{k}}D_k V_1 (w^H(m), w^H (m))_t
&=&
V_1 ( k (m), \hat{k} (m) )_t.
\end{eqnarray*}
Moreover, as $m \to\infty$, the right-hand sides of the above
equations converge~to
\begin{eqnarray*}
V_1 ( k , w^H )_t
\quad
\mbox{and}
\quad
V_1 ( k , \hat{k} )_t
\end{eqnarray*}
almost surely and in $L^2(\mu^H)$.
[Note that the above quantities are well-defined since $w^H$ is of
finite $p$-variation
and $k, \hat{k}$ is of finite $q$-variation with\vspace*{1pt} $1/p+1/q >1$.
Since $k, \hat{k}$ are of finite $(q-\varepsilon)$-variation for sufficiently
small $\varepsilon>0$, $k(m), \hat{k}(m)$ converge to $k, \hat{k}$
in $q$-variation norm, resp.]
\end{lm}

\begin{pf}
On $[ (l-1)/2^m, l/2^m ]$, $dw^H(m)_t = 2^n ( w^H_{ l/2^m } - w^H_{
(l-1)/2^m } ) \,dt$.
Therefore,
\[
D_k dw^H(m)_t = 2^n D_k \bigl( w^H_{ l/2^m } - w^H_{ (l-1)/2^m } \bigr) \,dt = 2^n
\bigl( k_{ l/2^m } - k_{ (l-1)/2^m } \bigr) \,dt
= dk(m)_t.
\]
From this, we see that
\begin{eqnarray*}
D_k \chi(w^H (m))_t
&=&
M_t \int_0^t M^{-1}_s \sigma( \phi^0_s) D_k \,dw^H(m)_s\\
&=&
M_t \int_0^t M^{-1}_s \sigma( \phi^0_s) \,dk(m)_s
=
\chi(k(m))_t.
\end{eqnarray*}
Since $\| k(m) -k \|_{q\mbox{-}\mathrm{var}}$ as $m\to\infty$ and
$k \mapsto\chi(k)$ is bounded linear\vspace*{1pt} from $C_0^{q\mbox{-}\mathrm{var} } ({\mathbf{R}}^d)$
to $C_0^{q\mbox{-}\mathrm{var} } ({\mathbf{R}}^d)$,
$\| \chi( k(m) )- \chi(k)\|_{q\mbox{-}\mathrm{var}}$ as $m \to\infty$.
[Note that, for sufficiently small $\varepsilon>0$, $k \in
C^{(q-\varepsilon)\mbox{-}\mathrm{var}}_0$
still holds.]
In a similar way,
\begin{eqnarray*}
&&\frac12
D_k V_1 (w^H (m), w^H(m) )_t\\
&&\qquad=
M_t \int_0^t M^{-1}_s
\{
\nabla\sigma(\phi^0_s ) \langle D_k \chi(w^H(m))_s , dw^H(m)_s
\rangle
\\
&&\hspace*{67pt}\qquad{}+\nabla\sigma(\phi^0_s ) \langle\chi( w^H(m))_s, D_k\, dw^H(m)_s
\rangle
\}
\\
&&\qquad=
M_t \int_0^t M^{-1}_s
\{
\nabla\sigma(\phi^0_s ) \langle\chi(k(m))_s , dw^H(m)_s \rangle
\\
&&\hspace*{65pt}\qquad{}+
\nabla\sigma(\phi^0_s ) \langle\chi( w^H(m))_s, d k(m)_s \rangle
\}
\\
&&\qquad= V_1(k(m), w^H (m))_t.
\end{eqnarray*}
Since $\| w^H (m) - w^H\|_{p\mbox{-}\mathrm{var} } \to0$ as $m\to\infty$ almost\vspace*{1pt} surely
and in $L^r $ for any $r>0$
(see~\cite{mss}),
$(1/2)D_k V_2 (w^H (m), w^H(m) )_t \to V_1 (k, w^H )_t$ almost surely
and in $L^2$.
Finally,
\[
(1/2)D_{\hat{k}}D_k V_2 (w^H (m), w^H(m) )_t
=
V_1 ( k (m), \hat{k} (m) )_t,
\]
which is nonrandom and clearly converges to $V_1 ( k , \hat{k} )_t$ as
$m \to\infty$.
\end{pf}

%pr5.2 #&#
\begin{pr}\label{prlimV2}
Let $V_1$ be as in (\ref{eqhessexp}) and consider $V_1 (w^H(m),w^H(m))^i_t$.
Here, $i$ stands for the $i$th component ($1 \le i \le n$).
Then, for each fixed $t$,
$V_1 (w^H(m),w^H(m))^i_t$ converges almost surely and in $L^2(\mu^H)$
as $m\to\infty$.
More precisely,
\[
\lim_{m \to\infty} V_1 (w^H(m),w^H(m))_t^i = \Theta^i_t
+
\Lambda^i_t.
\]
Here, $\Theta_t^i$ is an element in $\mathcal{C}_2$ which corresponds
to the symmetric Hilbert--Schmidt bilinear form $V_1 ( \bullet, \bullet)_t^i
$
and $t \mapsto\Lambda_t^i := \lim_{m \to\infty} E [ V_1
(w^H(m),w^H(m))_t^i]$ is of finite $p$-variation.
\end{pr}

\begin{pf}
First note that $V_1(x,x)$ has a rough path representation.
Recall the (stochastic) Taylor expansion of the It\^o map (\ref
{rdeqvar}) around $\gamma$.
Then, $V_1(x,x)$ was essentially calculated in computation for the
second Taylor term.
There is a continuous map
$V' \dvtx  G\Omega_p ({\mathbf{R}}^d) \to G\Omega_p ({\mathbf{R}}^n)$
such that\break
$V_1(x,x) = V' (X)^1$ for all $x \in C^{q\mbox{-}\mathrm{var}}_0 ({\mathbf{R}}^d)$.
Here, the superscript means the first level path and $X \in G\Omega_p
({\mathbf{R}}^d)$ is the lift of $x$.
Moreover, since the integral that defines $V_1$ or $V'$ in (\ref
{eqhessexp}) is of second order,
$V'$ has the following property: there exists a constant $c>0$ such
that, for all $X, Y \in G\Omega_p ({\mathbf{R}}^d)$,
\begin{eqnarray*}
\| V' (X)^1\|_{p\mbox{-}\mathrm{var}} &\le& c \bigl(1 + \xi(X)^2\bigr),
\\
\| V' (X)^1- V'(Y)^1\|_{p\mbox{-}\mathrm{var}} &\le& c \bigl(1 + \xi(X)^c\bigr) \sum_{j=1}^{[p]}
\| X^j-Y^j\|_{p/j\mbox{-}\mathrm{var}} .
\end{eqnarray*}
Here, $\xi(X) = \sum_{j=1}^{[p]} \| X^j\|_{p/j\mbox{-}\mathrm{var}}^{1/j}$.
From this, a.s.-convergence of $V_1 (w^H(m),\break w^H(m))=V'(w^H(m)))^1$ to
$V'(W^H))^1$ is obvious.

It is shown in~\cite{cq} that $E[ \| W^H(m)^j - W^{H,j}\|_{p/j\mbox{-}\mathrm{var}} ]
\to0$
as $m \to\infty$.
From Proposition~\ref{prfern},\vspace*{1pt} $\sup_m E[ \| W^H(m)^j \|_{p/j\mbox{-}\mathrm{var}}^r ]
<\infty$ for any $r>0$ and $1 \le j \le[p]$.
Then, we easily see from these and H\"older's inequality that
\[
E[ \| V'( W^H(m) )^1 - V' (W^{H} )^1 \|_{p\mbox{-}\mathrm{var}}^2 ] \to0 \qquad \mbox
{as $m \to\infty$.}
\]
This implies the $L^2$-convergence.
Since $ V' (W^{H} ) = \lim_{m \to\infty} V_1 (w^H(m),\break w^H(m))$ is a
$C^{p\mbox{-}\mathrm{var}}_0 ({\mathbf{R}}^n)$-valued
random variable,
\[
\|E[ V' (W^{H} )^1 ]\|_{p\mbox{-}\mathrm{var}}  \le E[ \| V' (W^{H} )^1 \|_{p\mbox{-}\mathrm{var}} ]
<\infty,
\]
which shows that $\Lambda$ is of finite $p$-variation.

By Lemma~\ref{lmclos} and the closability of the derivative operator
$D$ in $L^2(\mu^H)$,
\[
\tfrac12 D_k V' (W^H)_t^{1,i}
=
V_1 ( k , w^H )_t^i,
\qquad
\tfrac12 D_{\hat{k}} D_k V' (W^H)_t^{1,i}
=
V_1 ( k , \hat{k} )_t^i,
\]
where the superscript $i$ denotes the $i$th component of ${\mathbf{R}}^n$.
These equality imply that $V' (W^H)_t^{1,i} - E[ V' (W^H)_t^{1,i} ]$ is
in $\mathcal{C}_2$,
which corresponds to
$V_1 ( \bullet, \bullet)_t^i$.
\end{pf}

%le5.3 #&#
\begin{lm}
\label{lm2ndHS}
Let $p' >p$ and $F\dvtx  C_0^{p' \mbox{-}\mathrm{var}}( {\mathbf{R}}^n)$ be a Fr\'echet
differentiable function.
Let
\[
\Theta_t = \lim_{m \to\infty} V_1 (w^H(m),w^H(m))_t - E \Bigl[ \lim_{m
\to
\infty} V_1 (w^H(m),w^H(m))_t \Bigr]
\]
be as in Proposition~\ref{prlimV2}.
Then, $\nabla F (\phi^0) \langle\Theta\rangle\in\mathcal{C}_2 (\mu^H)$
which corresponds to the symmetric Hilbert--Schmidt bilinear form
$\nabla F (\phi^0) \circ V_1 = \nabla F (\phi^0) \langle V_1 (
\bullet,
\bullet) \rangle$ on~$\mathcal{H}^H$.
\end{lm}

\begin{pf}
Denote by $g_K$ the element of $\mathcal{C}_2 (\mu^H)$ which corresponds
to a symmetric
Hilbert--Schmidt bilinear form (or, equivalently, operator)
$K$ and
set
\[
M :=\{ \alpha\in C_0^{p \mbox{-}\mathrm{var}}( {\mathbf{R}}^n)^*  |  \alpha\langle\Theta
(w)\rangle= g_{\alpha\circ V_1 } (w)
 \mbox{ a.a. $w$  ($\mu^H$)} \}.
\]
Obviously, $M$ is a linear subspace.
Moreover, from Lemma~\ref{lmhilschV1}, $M$ is closed under weak*-limit.
By Lemma~\ref{prlimV2},
the evaluation map $\operatorname{ev}_t^i\  (t \in[0,1], 1 \le i \le n)$ defined by
$\operatorname{ev}_t^i \langle y \rangle= y^i_t$ is in $M$.
Denote by $\pi_m \dvtx  C_0^{p\mbox{-}\mathrm{var}}( {\mathbf{R}}^n) \to C_0^{p \mbox{-}\mathrm{var}}( {\mathbf{R}}^n)$ the projection defined by
$\pi(m) y =y(m)$, where $y(m)$ is the $m$th dyadic piecewise linear
approximation of $y \in C_0^{p\mbox{-}\mathrm{var}}( {\mathbf{R}}^n)$.
Note that $\nabla F(\phi^0) \langle\pi(m) y \rangle$ can be written
as a
linear combination of $y_{k/2^m}^i\ (1 \le k \le2^m, 1 \le i \le n)$.
Hence,\vspace*{-1pt} $\nabla F(\phi^0) \circ\pi(m) \in M$.
Since $p' >p$, $y(m) \to y$ in $p'$-variation norm.
This implies that $\nabla F(\phi^0) \circ\pi(m) \to\nabla F(\phi^0)$
in the weak*-topology.
Hence, $ \nabla F(\phi^0) \in M$.
\end{pf}

%%%%%%%%%%
Let $A_1$ be a self-adjoint Hilbert--Schmidt operator on $\mathcal{H}^H$
which corresponds to
\[
\nabla F (\phi^0) \langle V_1 ( \bullet, \bullet) \rangle.
\]
Then, $A- A_1$ is a self-adjoint Hilbert--Schmidt operator on $\mathcal{H}^H$ which corresponds to
\[
\nabla F (\phi^0) \langle V_2 ( \bullet, \bullet) \rangle+ \nabla
^2 F
(\phi
^0) \langle\chi(\bullet), \chi(\bullet) \rangle.
\]
Obviously, this bilinear form extends to one on $C_0^{p\mbox{-}\mathrm{var}} ({\mathbf{R}}^d)$ and, hence, is of trace class by Goodman's theorem.
See (\ref{eqhessexp}) for the definition\vadjust{\goodbreak} of $V_1, V_2$.
Combining these all, we see that
\[
k \in\mathcal{H}^H \mapsto\langle Ak,k\rangle_{\mathcal{H}^H}
=
\nabla F (\phi^0) \langle\psi( k, k) \rangle+ \nabla^2 F (\phi^0)
\langle
\chi
(k), \chi(k) \rangle
\]
extends to a continuous map on $G\Omega_p ({\mathbf{R}}^d)$ and we
denote it by $\langle AX,X\rangle$ for $X \in G\Omega_p ({\mathbf{R}}^d)$.
%

%
% Constant factor 1/2 was forgotten
%

%le5.4 #&#
\begin{lm}
\label{lmdet2}
Let $\alpha\ge1$ be such that $\mathrm{Id}_{\mathcal{H}^H} +\alpha A$ is
strictly positive in the form sense.
Then,
\begin{eqnarray*}
&&\int_{ \mathcal{X}} \exp\biggl( - \frac{\alpha}{2} \langle AW^H,W^H\rangle\biggr)
\mu^H (dw^H)
\\
&&\qquad=
\int_{ G\Omega_p ({\mathbf{R}}^d)} \exp\biggl( - \frac{\alpha}{2} \langle
AX,X\rangle\biggr)
{\mathbb P}^H (dX)
<\infty.
\end{eqnarray*}
In particular, $e^{ - \langle A \bullet, \bullet\rangle/2}$ is in
$L^r (
G\Omega_p ({\mathbf{R}}^d) , {\mathbb P}^H)$ for some $r>1$.
\end{lm}

\begin{pf}
As a functional of $w^H$,
$\langle(A-A_1)W^H,W^H \rangle$ is a sum of $\operatorname{Tr} (A-A_1)$
and the second order Wiener chaos corresponding to $A-A_1$.
From Proposition~\ref{prlimV2} and Lemma~\ref{lm2ndHS},
$\langle AW^H,W^H \rangle$ is a sum of a constant $\operatorname{Tr} (A-A_1) +
\nabla
F(\phi^0) \langle\Lambda\rangle$
and the second order Wiener chaos corresponding to $A$ (which is
denoted by $\Xi_A$ below).
It is well known (see Remark~\ref{recfdet} below) that
\[
{\mathbb E} [ e^{-\alpha\Xi_A /2}] = \det_2 ( \mathrm{Id}_{\mathcal{H}^H}
+\alpha A)^{-1/2},
\]
where $\det_2$ stands for the Carleman--Fredholm determinant.
\end{pf}

%re5.5 #&#
\begin{re}\label{recfdet}
Let $( \hat\mathcal{X} , \hat{H}, \hat\mu)$ be any abstract Wiener space.
For a symmetric Hilbert--Schmidt operator $\hat{A} \dvtx  \hat{H} \to\hat
{H}$, we denote by $\hat\Xi$ the corresponding
element in the second Wiener chaos $\hat\mathcal{C}_2$.
If $\hat{A} > - \mathrm{Id}$ in the form sense, then
%
%e5.1 #&#
\begin{equation}\label{eqredet2}
{\mathbb E} [ e^{- \hat\Xi/2}] = \det_2 ( \mathrm{Id} + \hat{A})^{-1/2}
\qquad
 \Biggl( :=
\prod_{j=1}^{\infty} \{ (1+ \lambda_j) e^{ - \lambda_j} \}^{-1/2}
\Biggr).
\end{equation}
Here, $\{\lambda_j \}_{j=1,2,\ldots}$ are eigenvalues of $\hat{A}$.

This fact is well known. For the reader's convenience, however, we give
a simple (and somewhat heuristic) proof below.
Suppose $ \hat\mathcal{X} = \hat{H} = {\mathbf{R}}^l$ and $\hat\mu$ is the
standard normal distribution.
Let $\hat{A}$ be such that
$
\langle\hat{A} \xi, \eta\rangle= \sum_{j=1}^l \lambda_j \xi_j
\eta_j ,\break
(\xi, \eta\in{\mathbf{R}}^l).
$
We assume that $\lambda_j >-1$ for all $j$.
In this case, $\hat\Xi$ is given by the following Hermite polynomial:
$
\hat\Xi(u) = \sum_{j=1}^l \lambda_j (u_j^2 -1)\ (u \in{\mathbf{R}}^l).$
This clearly corresponds to~$\hat{A}$ because
$(1/2)D_{\xi}D_{\eta} \hat\Xi(u) = \langle\hat{A} \xi, \eta
\rangle$.
(Recall that this is the way we identified symmetric Hilbert--Schmidt
operators with elements of second Wiener chaos
in the previous argument.)
A simple calculation shows that
\begin{eqnarray*}
{\mathbb E} [ e^{- \hat\Xi/2}] &=& \int_{ {\mathbf{R}}^l } e^{- \hat\Xi
(u) /2} \prod_{j=1}^{l} \frac{1}{\sqrt{2\pi} } e^{-u_j^2/2} \,du_j
\\
&=&
\prod_{j=1}^{l} \frac{1}{\sqrt{2\pi} } \int_{ {\mathbf{R}}^l } \exp\biggl(
\frac{
\lambda_j - (1+ \lambda_j ) u_j^2 }{2} \biggr) \,du_j
=
\prod_{j=1}^{l} \{ (1+ \lambda_j) e^{ - \lambda_j} \}^{-1/2}.
\end{eqnarray*}

We can easily do a similar computation in the case of $({\mathbf{R}}^{\infty
}, l^2, \mu^{\infty})$, where $\mu^{\infty}$ denotes the
countable product of the one-dimensional standard normal distribution.
The general case reduces to the case of ${\mathbf{R}}^{\infty}$, after
$\hat
{A}$ is diagonalized.

%%%%%%%%%%%%%%%%
Another method to verify (\ref{eqredet2}) is to use an explicit
formula for
the characteristic function of the quadratic Wiener functional, which
has been studied extensively.
Let $B \dvtx  \hat{H} \to\hat{H}$ be any symmetric Hilbert--Schmidt
operator and let $\Xi_B$ be the corresponding
element in the second Wiener chaos $\hat\mathcal{C}_2$.
Then, we have
\[
\int_{ \hat{\mathcal{X}} } \exp( (\zeta/2) \Xi_B ) \,d\hat
{\mu} =
\det_2 ( \mathrm{Id} -\zeta B)^{-1/2}
\qquad
\mbox{for any $\zeta\in\mathbf{C}$ with $|\zeta| < 1/\|B\|_{\mathrm{op}} $.}
\]
For example, see Janson~\cite{jan}, page 78, or Taniguchi~\cite{tan},
page 13.
The formula~(\ref{eqredet2}) immediately follows from this.
\end{re}

%%%%%%%%%%%%%%%%%%%%%%%%%%%%%%%%%%%%%%%%%%%%%%%%
%%%%%%%%%%%%%%%%%%%%%%%%%%%%%%%%%%%%%%%%%%%%%%
%% New Section
%%%%%%%%%%%%%%%%%%%%%%%%%%%%%%%%%%%%%%%%%%%%%%%%
%%%%%%%%%%%%%%%%%%%%%%%%%%%%%%%%%%%%%%%%%%%%%%%%%%%%%%%%%%
%s6 #&#
\section{Proof of Laplace approximation}\label{secproof}
%
%s6.1 #&#
\subsection{\texorpdfstring{Large deviation for the law of $Y^{\varepsilon}$ as $\varepsilon\searrow0$}
{Large deviation for the law of Y epsilon as epsilon -> 0}}\label{sec61}
In this section we prove the main theorem (Theorem~\ref{thmmain}).
Let $Y^{\varepsilon}$ be a solution of RDE (\ref{defrde}).
The law of $(Y^{\varepsilon})^1 = Y^{\varepsilon, 1}$ is the
probability measure on
$C_0^{p\mbox{-}\mathrm{var}} ({\mathbf{R}}^n)$ for any $p >1/H$.
Then, by Theorem~\ref{thmconti} and Proposition~\ref{prldp}
we can use the contraction principle to see that the law of $\{
Y^{\varepsilon
,1} \}_{\varepsilon>0}$
satisfies large deviation as $\varepsilon\searrow0$.
The good rate function is given as follows:
\[
I(y) =
\cases{
\inf\{ \| k \|_{\mathcal{H}^H}^2 /2  |  y =\hat\Phi_0 (k, \lambda)^1 \}
&\quad
$\mbox{(if $y =\hat\Phi_0 (k, \lambda)^1$ for some $k \in\mathcal{H}^H$),}$ \vspace*{2pt}\cr
\infty&\quad $\mbox{(otherwise).}$ }
\]
Here, $\hat\Phi_{\varepsilon}$ is the It\^o map corresponding to RDE
(\ref
{defrde}) and $\lambda_t =t$.
For a bounded continuous function $F$ on $C_0^{p\mbox{-}\mathrm{var}} ({\mathbf{R}}^n)$, it
holds that
\[
\lim_{\varepsilon\searrow0} \varepsilon^2 \log{\mathbb E}\bigl [ \exp
\bigl(- F(Y^{\varepsilon,1} )
/\varepsilon^2\bigr) \bigr]
=
-\inf\{ F(y) +I(y)  |  y \in C_0^{p\mbox{-}\mathrm{var}} ({\mathbf{R}}^n) \}.
\]

Now, let us consider Laplace's method, that is, the precise asymptotic
behavior of
the following integral:
\begin{eqnarray*}
{\mathbb E} \bigl[ \exp\bigl(- F(Y^{\varepsilon,1} ) /\varepsilon^2\bigr) \bigr]
&=&
\int_{ G\Omega_p ({\mathbf{R}}^d) }
\exp\bigl(- F( \hat\Phi_{\varepsilon} ( \varepsilon X ,\lambda)^1 )
/\varepsilon^2\bigr)
{\mathbb P}^{H} (dX)
\nonumber\\
&=&
\int_{ G\Omega_p ({\mathbf{R}}^d) }
\exp\bigl(- F( \hat\Phi_{\varepsilon} ( X ,\lambda)^1 ) /\varepsilon^2\bigr)
{\mathbb P}^{H}_{\varepsilon} (dX)
\nonumber
\end{eqnarray*}
as $\varepsilon\searrow0$ under Assumptions (H1)--(H4).

Let $\gamma\in\mathcal{H}^H \subset G\Omega_p ({\mathbf{R}}^d)$
be the unique element at which $F( \hat\Phi_0 (  \cdot  , \lambda
) )
+ \|  \cdot \|^2_{\mathcal{H}^H} /2$
attains minimum ($F_{\Lambda} (\gamma)=: a$) as in (H2).
By a well-known argument, for any neighborhood of $O \subset G\Omega_p
({\mathbf{R}}^d)$ of $\gamma$,
there exist positive constants $\delta, C$ such that
\[
\int_{ O^c }
\exp\bigl(- F( \hat\Phi_{\varepsilon} ( X ,\lambda)^1 ) /\varepsilon^2\bigr)
{\mathbb P}^{H}_{\varepsilon} (dX)
\le
C e^{-(a+\delta) /\varepsilon^2},
\qquad
\varepsilon\in(0,1].
\]
This decays very fast and does not contribute to the asymptotic expansion.

%%%%%%%%%%%%%%%%%%%%%%%%%%%%%%%%%%%%%%%%%%%%%
%s6.2 #&#
\subsection{\texorpdfstring{Computation of $\alpha_0$}{Computation of alpha 0}}\label{subsec2nd}
In this subsection we compute the first term $\alpha_0$ in the
asymptotic expansion when $G \equiv1$ (constant)
and show $\alpha_0 >0$.
To do so, we need the (stochastic) Taylor expansion (Theorem \ref
{thmtaylor}) up to order $m=2$.
Once this is done,
expansion up to higher order terms can be obtained rather easily.

For $\rho>0$, set $U_{\rho}=\{ X \in G\Omega_p ({\mathbf{R}}^d)  |  \xi(X)
< \rho\}$,
where $\xi$ is given in (\ref{defxi}).
Then, taking $O = \gamma+ U_{\rho}$,
we see from the theorem of the Cameron--Martin type (Proposition \ref
{prCM}) that
%
%e6.1 #&#
\begin{eqnarray}\label{eqlocint}
&&
\int_{ \gamma+U_{\rho} }
\exp\bigl(- F( \hat\Phi_{\varepsilon} ( X ,\lambda)^1 ) /\varepsilon^2\bigr)
{\mathbb P}^{H}_{\varepsilon} (dX)
\nonumber\\
&&\qquad=
\int_{ U_{\rho} }
\exp\bigl(- F\bigl( \hat\Phi_{\varepsilon} ( X +\gamma,\lambda)^1 \bigr)
/\varepsilon^2\bigr)
\nonumber
\\[-8pt]
\\[-8pt]
\nonumber
&&\hspace*{28pt}\qquad{}\times \exp\biggl( -\frac{1}{\varepsilon^2} \langle\gamma, X^1 \rangle-\frac
{1}{2\varepsilon^2} \|
\gamma\|^2_{\mathcal{H}^H} \biggr)
{\mathbb P}^{H}_{\varepsilon} (dX)
\\
&&\qquad=
\int_{ \{\xi(\varepsilon X) < \rho\} }
\exp\biggl(
- \frac{F( \phi^{(\varepsilon)} ) }{\varepsilon^2}
-\frac{1}{\varepsilon} \langle\gamma, X^1 \rangle-\frac
{1}{2\varepsilon^2} \| \gamma\|
^2_{\mathcal{H}^H} \biggr)
{\mathbb P}^{H} (dX).\nonumber
\end{eqnarray}
As we will see, $ \langle\gamma,  \cdot  \rangle$ extends to a continuous
linear functional on $C_0^{p\mbox{-}\mathrm{var}} ({\mathbf{R}}^d)$
and, in particular, everywhere defined.

For sufficiently small $\rho$ (i.e., $\rho\le\rho_0$ for some $\rho
_0$), $ \phi^{(\varepsilon)} $ is in the neighborhood of $\phi^0$ as in
Assumption (H3).
So, from the Taylor expansion for~$F$,
\begin{eqnarray}
F( \phi^{(\varepsilon)} )
&=&
F(\phi^0) + \nabla F (\phi^0) \bigl\langle\phi^{(\varepsilon)} - \phi
^0 \bigr\rangle
+ \frac12 \nabla^2 F (\phi^0) \bigl\langle\phi^{(\varepsilon)} - \phi
^0, \phi
^{(\varepsilon)}
- \phi^0 \bigr\rangle
\nonumber\\
&&{}
+\frac16 \int_0^1 \,d\theta\nabla^3 F\bigl ( \theta\phi^{(\varepsilon)} +(1-
\theta
) \phi^0\bigr)
\bigl\langle\phi^{(\varepsilon)} - \phi^0, \phi^{(\varepsilon)} - \phi
^0, \phi^{(\varepsilon)}
- \phi
^0 \bigr\rangle
\nonumber\\
&=&
F(\phi^0) + \nabla F (\phi^0) \langle\varepsilon\phi^1 +
\varepsilon^2 \phi^2 \rangle
+ \frac12 \nabla^2 F (\phi^0) \langle\varepsilon\phi^1,
\varepsilon\phi^1 \rangle+
Q^3_{\varepsilon}.
\nonumber
\end{eqnarray}
Here, the remainder term $Q^3_{\varepsilon}$ satisfies the following estimates:
there exists a positive constant $C=C(\rho_0)$ such that
%
%e6.2 #&#
\begin{equation} \label{estQ3}
| Q^3_{\varepsilon} | \le C \bigl(\varepsilon+\xi(\varepsilon X)\bigr)^3
\qquad  \mbox{on the set
$\{\xi(\varepsilon X) < \rho_0 \} $. }
\end{equation}
Note that $C$ is independent of the choice of $\rho\ (\rho\le\rho_0)$.

Now we compute the shoulder of $\exp$ on the right-hand side of (\ref
{eqlocint}).
Terms of order $-2$ are computed as follows:
\[
- \frac{1}{\varepsilon^2} \biggl( F( \phi^0 ) +\frac12 \| \gamma\|
^2_{\mathcal{H}^H}\biggr)
=- \frac{a}{\varepsilon^2} .
\]

Since $k \in\mathcal{H}^H \mapsto F( \Phi_0 (k, \lambda) ) +\| k\|
^2_{\mathcal{H}^H} /2$ takes its minimum at $k=\gamma$,
we see that
\[
\langle k, \gamma\rangle_{\mathcal{H}^H} + \nabla F( \phi^0) \langle
\chi(k) \rangle=0,
\]
where $ \chi(k)$ is given by (\ref{defchik}) or (\ref{defchik2}).
By (\ref{defchik2}) and the Young integral,
$k \mapsto\nabla F( \phi^0) \langle\chi(k) \rangle$ extends to a continuous
linear map from $C_0^{p\mbox{-}\mathrm{var}} ({\mathbf{R}}^d)$
and so does $\langle\gamma,  \cdot \rangle_{\mathcal{H}^H} $.
Hence, the measurable linear functional (i.e., the first Wiener chaos)
associated with $\gamma$
is this continuous extension.

An ODE for $\phi^1 =\phi^1 (k)=\phi^1 (k, \gamma)$ is as follows [$k
\in C_0^{q\mbox{-}\mathrm{var}} ({\mathbf{R}}^d)$]:
%
%e6.3 #&#
\begin{eqnarray}\label{defphi1}
&&d \phi^1_t - \nabla\sigma(\phi^0_t) \langle\phi^1_t, d\gamma_t
\rangle-
\nabla_y \beta(0 ,\phi^0_t) \langle\phi^1_t \rangle \,dt
\nonumber
\\[-8pt]
\\[-8pt]
\nonumber
&&\qquad=
\sigma(\phi^0_t) \,dk_t + \nabla_{\varepsilon} \beta(0, \phi^0_t) \,dt,
 \qquad \phi^1_0=0.
\end{eqnarray}
Note that both $\phi^1$ and $\chi$ extend to a continuous map from
$G\Omega_p ({\mathbf{R}}^d)$.
The difference $\theta^1_t := \phi^1_t (X)-\chi_t (X)$ is independent
of $X$ (i.e., nonrandom), of finite $q$-variation, and satisfies
%
%e6.4 #&#
\begin{eqnarray}\label{deftheta1}
&&d \theta^1_t - \nabla\sigma(\phi^0_t) \langle\theta^1_t, d\gamma_t
\rangle
- \nabla_y \beta(0 ,\phi^0_t) \langle\theta^1_t \rangle \,dt
\nonumber
\\[-8pt]
\\[-8pt]
\nonumber
&&\qquad=
\nabla_{\varepsilon} \beta(0, \phi^0_t) \,dt,\qquad
 \theta^1_0=0.
\end{eqnarray}
Or, equivalently, $ \theta^1_t = M_t \int_0^t M^{-1}_s \nabla
_{\varepsilon}
\beta(0, \phi^0_s) \,ds$.
Consequently, terms of order $-1$ are computed as follows:
\[
-\frac{1}{\varepsilon} \bigl(
\nabla F (\phi^0) \langle\phi^1\rangle+ \langle\gamma, X^1\rangle
\bigr)
=
-\frac{\nabla F (\phi^0) \langle\theta^1\rangle}{\varepsilon}.
\]

Now we compute terms of order $0$.
The second term $\phi^2 = \phi^2 (k)=\phi^2 (k, \gamma)$ in the
expansion in Theorem~\ref{thmtaylor} satisfies the following ODE (see
\cite{ina2}, e.g.):
%
%e6.5 #&#
\begin{eqnarray}\label{eqodephi2}
&&
d \phi^2_t - \nabla\sigma(\phi^0_t) \langle\phi^2_t, d\gamma_t
\rangle-
\nabla_y \beta(0 ,\phi^0_t) \langle\phi^2_t \rangle \,dt
\nonumber\\
&&\qquad=
\nabla\sigma(\phi^0_t) \langle\phi^1_t, dk_t \rangle
+
\tfrac12 \nabla^2 \sigma(\phi^0_t) \langle\phi^1_t, \phi
^1_t,d\gamma
_t \rangle
\nonumber
\\[-8pt]
\\[-8pt]
\nonumber
&&\quad\qquad{}+
\tfrac12 \nabla_y^2 \beta(\phi^0_t) \langle\phi^1_t, \phi
^1_t\rangle \,dt
\\
&&\qquad\quad{}
+
\nabla_y \nabla_{\varepsilon} \beta(\phi^0_t) \langle\phi^1_t
\rangle \,dt
+\tfrac12 \nabla_{\varepsilon}^2 \beta(0, \phi^0_t) \,dt,\qquad
 \phi^2_0=0.\nonumber
\end{eqnarray}

Let $\chi$ and $\psi$ be as in (\ref{defchik}) and (\ref{defpsi}),
respectively.
By the same argument for (stochastic) Taylor expansion (Theorem \ref
{thmtaylor}),
those extend to continuous maps from $G\Omega_p ({\mathbf{R}}^d)$ and we
write $\chi(X)$ and $\psi(X)=\psi(X,X)$.
If we set $\theta^2 (k):= \phi^2(k) - \psi(k) /2$, then $\theta^2$
satisfies the following ODE:
%
%e6.6 #&#
\begin{eqnarray}\label{deftheta2}
&&
d \theta^2_t - \nabla\sigma(\phi^0_t) \langle\theta^2_t, d\gamma_t
\rangle
- \nabla_y \beta(0 ,\phi^0_t) \langle\theta^2_t \rangle \,dt
\nonumber\\
&&\qquad=
\nabla\sigma(\phi^0_t) \langle\theta^1_t, dk_t \rangle
+
\tfrac12 \nabla^2 \sigma(\phi^0_t) \langle\theta^1_t, \theta
^1_t,d\gamma
_t \rangle
+
\nabla^2 \sigma(\phi^0_t) \langle\theta^1_t, \chi_t,d\gamma_t
\rangle
\nonumber
\\[-8pt]
\\[-8pt]
\nonumber
&&\qquad\quad{}
+
\tfrac12 \nabla_y^2 \beta(\phi^0_t) \langle\theta^1_t, \theta
^1_t\rangle \,dt
+
\nabla_y^2 \beta(\phi^0_t) \langle\theta^1_t, \chi_t\rangle \,dt
\\
&&\qquad\quad{}
+
\nabla_y \nabla_{\varepsilon} \beta(\phi^0_t) \langle\theta^1_t
+\chi_t\rangle \,dt
+\tfrac12 \nabla_{\varepsilon}^2 \beta(0, \phi^0_t) \,dt,\qquad
 \theta^2_0=0.\nonumber
\end{eqnarray}
Or, equivalently,
\begin{eqnarray*}
\theta^2_t
&=&
M_t \int_0^t M^{-1}_s \biggl(
\nabla\sigma(\phi^0_s) \langle\theta^1_s, dk_s \rangle
+
\frac12 \nabla^2 \sigma(\phi^0_s) \langle\theta^1_s, \theta
^1_t,d\gamma
_s \rangle
\\
&&\hspace*{36pt}\qquad{}+
\nabla^2 \sigma(\phi^0_s) \langle\theta^1_s, \chi_s , d\gamma_s
\rangle
\\
&&\hspace*{58pt}{}+
\frac12 \nabla_y^2 \beta(\phi^0_s) \langle\theta^1_s , \theta^1_s
\rangle \,ds
+
\nabla_y^2 \beta(\phi^0_s) \langle\theta^1_s , \chi_t\rangle \,ds
\\
&&\hspace*{61pt}{}+
\nabla_y \nabla_{\varepsilon} \beta(\phi^0_s ) \langle\theta^1_t
+\chi_s
\rangle \,ds
+
\frac12 \nabla_{\varepsilon}^2 \beta(0, \phi^0_s) \,ds \biggr).
\end{eqnarray*}
This is just a Young integral and $k \mapsto\theta^2(k)$ extends to a
continuous map from $C_0^{p\mbox{-}\mathrm{var}} ({\mathbf{R}}^d)$
or from $G\Omega_p ({\mathbf{R}}^d)$ to $C_0^{p\mbox{-}\mathrm{var}} ({\mathbf{R}}^n)$.
Moreover, $\theta^2$ is of first order, that is, for some constant $C>0$,
$\| \theta^2 (X) \|_{p\mbox{-}\mathrm{var}} \le C (1 + \xi(X))$ holds for any $X \in
G\Omega_p ({\mathbf{R}}^d)$.
In particular, by the Fernique-type theorem (Proposition~\ref{prfern}),
(a~constant multiple of) $\theta^2$ is exponentially
integrable.

Hence, terms of order $0$ on the shoulder of $\exp$ on the right-hand
side of (\ref{eqlocint})
are as follows:
%
%e6.7 #&#
\begin{eqnarray}\label{eqzeroth}
&&\nabla F (\phi^0) \langle\phi^2 \rangle+ \tfrac12 \nabla^2 F(\phi
^0) \langle
\phi
^1, \phi^1 \rangle\nonumber\\
&&\qquad=
\tfrac12 [
\nabla F(\phi^0) \langle\psi\rangle+ \nabla^2 F(\phi^0) \langle
\chi, \chi
\rangle]
%&&
+
\nabla F(\phi^0) \langle\theta^2 \rangle\\
&&\qquad\quad{}
+\tfrac12 \nabla^2 F(\phi^0) \langle\theta^1, \theta^1 \rangle
+ \nabla^2 F(\phi^0) \langle\theta^1, \chi\rangle.\nonumber
\end{eqnarray}
Note that the last three terms on the right-hand side are dominated by
$C (1 +\xi(X))$
and that the first term is $\langle AX,X\rangle/2$ as in Lemma~\ref{lmdet2}.
By Proposition~\ref{prfern} and Lemma~\ref{lmdet2},
\begin{eqnarray}
\exp\bigl(
- \nabla F (\phi^0) \langle\phi^2 \rangle- \tfrac12 \nabla^2
F(\phi^0)
\langle\phi
^1, \phi^1 \rangle\bigr)
\in
L^r ( G\Omega_p ({\mathbf{R}}^d) , {\mathbb P}^H) \nonumber\\
\eqntext{\mbox{for some $r>1$.}}
\end{eqnarray}
If $\rho>0$ is chosen sufficiently small, then $\exp[ 2C\rho(1 +\xi
(X) )^2 ] \in L^{r'} ( G\Omega_p ({\mathbf{R}}^d) , \break {\mathbb P}^H)$
for the conjugate exponent $r'$, that is, $1/r + 1/r' =1$. (We
determine $\rho$, here.)
We easily see that, if $\varepsilon\le\rho$,
%
%e6.8 #&#
\begin{eqnarray}
\qquad &&\mathbf{1}_{ \{ \xi(\varepsilon X) <\rho\} }
\exp\bigl(
- \nabla F (\phi^0) \langle\phi^2 \rangle- \tfrac12 \nabla^2
F(\phi^0)
\langle\phi
^1, \phi^1 \rangle\bigr)
\exp( - \varepsilon^{-2} Q^3_{\varepsilon} )
\nonumber
\\[-8pt]
\\[-8pt]
\nonumber
&&\qquad\le
\exp\bigl(
- \nabla F (\phi^0) \langle\phi^2 \rangle- \tfrac12 \nabla^2
F(\phi^0)
\langle\phi
^1, \phi^1 \rangle\bigr)
\exp\bigl[ 2C\rho\bigl(1 +\xi(X) \bigr)^2 \bigr].
\end{eqnarray}
The right-hand side is integrable and independent of $\varepsilon$.
So, we may use the dominated convergence theorem to obtain that
\begin{eqnarray*}
&&\lim_{ \varepsilon\searrow0}
\int_{ \{ \xi(\varepsilon X) <\rho\} }
\exp\biggl(
- \nabla F (\phi^0) \langle\phi^2 \rangle- \frac12 \nabla^2
F(\phi^0)
\langle\phi
^1, \phi^1 \rangle- \frac{1}{ \varepsilon^{2}} Q^3_{\varepsilon}
\biggr)
{\mathbb P}^H (dX)
\\
&&\qquad=
\int_{ G\Omega_p ({\mathbf{R}}^d) }
\exp\biggl(
- \nabla F (\phi^0) \langle\phi^2 \rangle- \frac12 \nabla^2
F(\phi^0)
\langle\phi
^1, \phi^1 \rangle\biggr)
{\mathbb P}^H (dX).
\end{eqnarray*}
By Lemma~\ref{lmdet2}, the right-hand side exists.
Thus, we have computed (the asymptotics of) (\ref{eqlocint}) up to
$\alpha_0$.

%%%%%%%%%%%%%%%%%%%%%%%%%%%%%%%%%%%%%%%%%%%%
%s6.3 #&#
\subsection{Asymptotic expansion up to any order}\label{subsecgente}
In this subsection we obtain the Laplace asymptotic expansion up to any order.
Since this is routine once $\alpha_0$ is obtained, we only give a
sketch of the proof.

By combining the (stochastic) Taylor expansions for $F, G$ and $\phi
^{(\varepsilon)}$,
we get
\begin{eqnarray*}
F\bigl( \phi^{(\varepsilon)} \bigr) - F( \phi^0) &\sim& \varepsilon\eta^1
+\cdots+ \varepsilon^n
\eta^n
+Q_{\varepsilon}^{n+1}\qquad \mbox{as $\varepsilon\searrow0$},
\\
G\bigl( \phi^{(\varepsilon)} \bigr) - G( \phi^0) &\sim& \varepsilon\hat\eta
^1 +\cdots+\varepsilon
^n \hat
\eta^n + \hat{Q}_{\varepsilon}^{n+1}
\qquad\mbox{as $\varepsilon\searrow0$}.
\end{eqnarray*}
Here, the remainder terms $ Q_{\varepsilon}^{n+1}, \hat
{Q}_{\varepsilon}^{n+1}$ satisfy
similar estimates to (\ref{estQ3}).

From this we see that
%
%e6.9 #&#
\begin{eqnarray}
&&
\int_{ \gamma+U_{\rho} }
G ( \hat\Phi_{\varepsilon} ( X ,\lambda)^1 )
\exp\bigl(- F( \hat\Phi_{\varepsilon} ( X ,\lambda)^1 ) /\varepsilon^2\bigr)
{\mathbb P}^{H}_{\varepsilon} (dX)
\nonumber\\
&&\qquad=
e^{-a/\varepsilon^2} e^{- \nabla F(\phi^0) \langle\theta^1\rangle
/\varepsilon}
%
%&&
% \times
\nonumber
\\[-8pt]
\\[-8pt]
\nonumber
&&\qquad\quad{}
 \times\int_{ \{\xi(\varepsilon X) < \rho\} }
G\bigl( \phi^{(\varepsilon)}\bigr)
\exp\biggl(
-\nabla F(\phi^0) \langle\phi^2\rangle-\frac12 \nabla^2 F(\phi
^0) \langle
\phi
^1, \phi^1\rangle
\biggr)
\\
&&\hspace*{68pt}\qquad{}\times \exp(-Q^3_{\varepsilon} / \varepsilon^2)
{\mathbb P}^{H} (dX)\nonumber
\end{eqnarray}
can easily be expanded.
Note that
\[
\biggl| e^u - \biggl(1 +\frac{u}{1!} + \cdots+ \frac{u^{n-1}}{(n-1)!} \biggr)
\biggr| \le\frac{ e^{|u|} |u|^n}{n!}
\qquad
\mbox{(with $u = - Q^3_{\varepsilon} / \varepsilon^2$)}
\]
and that $Q^3_{\varepsilon}= \varepsilon^3 \eta^3 +\cdots+
\varepsilon^n \eta^n +Q_{\varepsilon}^{n+1}$.
Thus, we have shown the main theorem (Theorem~\ref{thmmain}).

%%%%%%%%%%%%%%%%%%%%%%%%%%%%%%%%%%%%%%%%%%%%%
%%%%%%%%%%%%%%%%%%%%%%%%%%%%%%%%%%%%%%%%%%%%%
%%%%%%%%%%%%%%%%%%%%%%%%%%%%%%%%%%%%%%%%%%
%s7 #&#
\section{Fractional order case: with an application to short time expansion}\label{sec7}
In this section we consider an RDE, which involves a fractional order
term of $\varepsilon$.
As a result, a\vadjust{\goodbreak} fractional order term of $\varepsilon$ appears in
the asymptotic expansion.
By time change, this has an application to the short time problems for
the solutions of the RDE
driven by fBRP.

%%%%%%%%%%%%%%%

First we see the scale invariance of fBRP.
It is well known that, for $0< c \le1$,
$( c^{-H} w^H_{c t} )_{0 \le t \le1}$ and $w^H$ have the same law.
A similar fact holds\vspace*{1pt} for the law of fBRP $W^H= (W^H_{s,t} )_{0 \le s
\le t \le1}$.
This is not so obvious from the
scale invariance of fBM $w^H$,
since fBRP $W^H$ is constructed via the dyadic partition of
$[0,1]$.\vspace*{-2pt}
%

%pr7.1 #&#
\begin{pr} \label{prscale}
Let $H \in(1/4 , 1/2)$ and $0< c\le1$.
Then, $( c^{-H}\times  W^H_{c s, c t} )_{0 \le s \le t \le1}$ and $W^H$ have
the same law.\vspace*{-2pt}
\end{pr}

\begin{pf}
(i)
Baudoin and Coutin showed this statement in~\cite{bc2}.\vspace*{-6pt}

\begin{longlist}
\item[(ii)]  Friz and Victoir~\cite{fvgauss} showed the following:
If a sequence of partitions of $[0,1]$ whose mesh tending to zero
satisfies a condition called ``nested,''
then the lift of $w^H$ via this sequence gives the same $W^H$ again.
Combining this result with the scaling property of $w^H$, we can
easily see the Proposition holds
at least for $c \in\mathbf{Q}$.
For $c \notin\mathbf{Q}$, just take a limit.\quad\qed\vspace*{-2pt}
\end{longlist}
\noqed\end{pf}

%%%%%%%%%%%

Let $H \in(1/4,1/3) \cup(1/3, 1/2)$.
For simplicity, we consider the following RDE:
%
%e7.1 #&#
\begin{equation}\label{eqrdefrac}
dY^{\varepsilon}_t = \sigma(Y^{\varepsilon}_t ) \varepsilon \,dX_t +
\varepsilon^{1/H} \hat\beta
(Y^{\varepsilon
}_t) \,dt,
\qquad
Y^{\varepsilon}_0 = 0.
\end{equation}
Here, $\sigma$ is as in Theorem~\ref{thmmain},
but we assume that a $C^{\infty}_b$-function $\hat\beta\dvtx {\mathbf{R}}^n
\to
{\mathbf{R}}^n$
and the drift term
is of this special form in this case.
Set $\beta(\varepsilon,y )= \varepsilon^{1/H} \hat\beta(y)$.
We also consider the following RDE, which is independent of
$\varepsilon$:
%
%e7.2 #&#
\begin{equation}\label{eqrdefrac222}
dV_t = \sigma(V_t ) \,dX_t + \hat\beta(V_t) \,dt,
\qquad
V_0 = 0.
\end{equation}
Basically, when we introduce randomness,
we always set $X= W^H$ in (\ref{eqrdefrac}) and (\ref{eqrdefrac222}).
Then, by the scale invariance of $W^H$ (see Proposition~\ref{prscale} below),
$(V_{ \varepsilon^{1/H}s, \varepsilon^{1/H} t})_{0 \le s \le t \le1}
$ and $(Y^{\varepsilon
}_{s,t})_{0 \le s \le t \le1} $ have the same law.
In particular, for each fixed $T \in(0,1]$, the ${\mathbf{R}}^n$-valued
random variables
$V^1_{ 0, T}$ and $(Y^{T^H})^1_{0,1}$
have the same law.
Therefore, the short time asymptotics for $V^1_{0,t}$ is related to the
small asymptotics of $(Y^{\varepsilon})^1$.

%%%%%%%%%%%%%%%%%%%%%
Let us fix some notation for fractional order expansions.
For
\[
M=\biggl\{ n_1 + \frac{n_2 }{H}  \Big|  n_1, n_2 = 0,1,2, \ldots\biggr\},
\]
let
$0=\kappa_0 < \kappa_1 <\kappa_2 <\cdots$
be all elements of $M$ in increasing order.
More concretely, leading terms are as follows:
%
%e7.3 #&#
\begin{eqnarray}\label{eqindex}
&&(\kappa_0, \kappa_1, \kappa_2, \ldots)
\nonumber\\
&&\qquad=
\biggl(0, 1,2, \frac{1}{H}, 3, 1+ \frac{1}{H} , 4, 2+ \frac{1}{H} ,
5\wedge
\frac{2}{H}, \ldots\biggr)\qquad
\mbox{if $H \in(1/3,1/2)$,}
\nonumber
\\[-8pt]
\\[-8pt]
\nonumber
\quad&&( \kappa_0, \kappa_1, \kappa_2, \ldots)
\\
&&\qquad=
\biggl(0, 1,2,3, \frac{1}{H}, 4, 1+ \frac{1}{H}, 5, \ldots, \biggr)\qquad
\mbox{if $H \in(1/4,1/3)$.}\nonumber
\end{eqnarray}
%
%It is important to note that, in (\ref{eqindex}), terms up to degree
%two (i.e., $\kappa_0, \kappa_1, \kappa_2$) are
%the same as in the previous sections.

As in the previous sections, we write $Y^{\varepsilon} = \hat\Phi
_{\varepsilon}(
\varepsilon
X)$, $\tilde{Y}^{\varepsilon} = \hat\Phi_{\varepsilon}(\varepsilon
X +\gamma)$,
and $\phi^{\varepsilon} =(\tilde{Y}^{\varepsilon} )^1$
for the solution of (\ref{eqrdefrac}).
By slightly modifying Theorem~\ref{thmtaylor}, we can prove the
(stochastic) Taylor expansion (around $\gamma$)
for
\[
\phi^{(\varepsilon)} = \phi^{0} +\varepsilon^{\kappa_1} \phi
^{\kappa_1} + \varepsilon
^{\kappa
_2} \phi^{\kappa_2} +\cdots+\varepsilon^{\kappa_m} \phi^{\kappa
_m} +
R^{\kappa
_{m+1}}_{\varepsilon}.
\]
In this case, $\phi^0$ satisfies the following ODE (in $q$-variation sense):
%
%e7.4 #&#
\begin{equation}\label{defqodefrac}
d\phi^0_t = \sigma( \phi^0_t) \,d\gamma_t,
\qquad
\phi^0_0 =0.
\end{equation}

%re7.2 #&#
\begin{re}
Although $(d/d\varepsilon)^m \vert_{\varepsilon=0}$ does not operate
on the
right-hand side of the following (formal) ODE,
%
%e7.5 #&#
\begin{equation}\label{eqrdefrac2}
d \phi^{(\varepsilon)}_t = \sigma\bigl(\phi^{(\varepsilon)} _t \bigr)
d(\varepsilon X_t +\gamma) +
\varepsilon
^{1/H} \hat\beta\bigl( \phi^{(\varepsilon)} _t\bigr) \,dt,\qquad
\tilde{Y}^{\varepsilon}_0 = 0,
\end{equation}
the proof of expansion in~\cite{ina2}, which is similar to Azencott's
argument in~\cite{az},
does not use the $\varepsilon$-derivative
and can be easily modified to our case.

Roughly and formally speaking, the proof goes as follows.
First, combine
\[
\phi^{(\varepsilon)} - \phi^{0} =\varepsilon^{\kappa_1} \phi
^{\kappa_1} +\cdots
+\varepsilon
^{\kappa_m} \phi^{\kappa_m} +\cdots
\]
and the Taylor expansion of $\sigma$ and $\hat\beta$ around $\phi^0_t$.
Next, pick up the terms of order $\alpha_m (m=1,2,\ldots)$.
Then, we obtain a very simple ODE of first order for $\phi^{\kappa_m}$
recursively.
This, in turn, can be used to
rigorously define $\phi^{\kappa_m}$.
In the end, we prove growth of the remainder term is of an expected order.
(This part is nontrivial and requires much computation.)
Note that this method can be used both in integer order and in
fractional order cases.
\end{re}

In the same way as in the previous sections, we have the following
modification of the main theorem
(Theorem~\ref{thmmain}).
%
%th7.3 #&#
\begin{tm}
\label{thmmainfrac}
Let the coefficients $\sigma\dvtx  {\mathbf{R}}^n \to\operatorname{Mat} (n,d)$ and $\hat
\beta\dvtx  {\mathbf{R}}^n \to{\mathbf{R}}^n$
be $C_b^{\infty}$
and consider the RDE (\ref{eqrdefrac}) with $X= W^H$, where $H \in
(1/4,1/3) \cup(1/3,1/2)$.
For simplicity, assume \textup{(H1)--(H4)} for any order $m$.
Then, we have the following asymptotic expansion as $\varepsilon
\searrow0$:
there are real constants $c$ and $\alpha_{\kappa_0} (= \alpha_0),
\alpha
_{\kappa_1}, \alpha_{\kappa_2},\ldots$ such that
\begin{eqnarray*}
&&
{\mathbb E} \bigl[
G( Y^{\varepsilon, 1}) \exp\bigl( - F ( Y^{\varepsilon, 1})
/\varepsilon^2 \bigr)
\bigr]
\\
&&\quad=
\exp\bigl( - F_{\Lambda} (\gamma) /\varepsilon^2 \bigr) \exp(- c/\varepsilon)
\cdot
\bigl(
\alpha_{\kappa_0} +\alpha_{\kappa_1} \varepsilon^{\kappa_1} +
\cdots+
\alpha
_{\kappa_m} \varepsilon^{\kappa_m} +O(\varepsilon^{\kappa_{m+1}})
\bigr)
\end{eqnarray*}
for any $m \ge0$.
\end{tm}

%re7.4 #&#
\begin{re}
It is important to note that, in (\ref{eqindex}), indices up to degree
two (i.e., $\kappa_0, \kappa_1, \kappa_2$) are
the same as in the previous sections.
The most difficult part of the proof of Theorem~\ref{thmmain}
is obtaining $\alpha_0$ [or checking that $\alpha_0 \in(0, \infty)$
when $G \equiv1$],
in which the (stochastic) Taylor expansion of $\phi^{(\varepsilon)}$
up to
$\phi
^2$ is used
(see Section~\ref{subsec2nd}).
Therefore, the proof in Section~\ref{subsec2nd} holds true without
modification in this case, too.
Higher order terms are different in the fractional order case. But, the
argument in Section~\ref{subsecgente}
is simple anyway and can easily be modified.
Thus, we can prove Theorem~\ref{thmmainfrac} without much difficulty.
\end{re}

%%%%%%%%%%%%%%%%%%

As a corollary, we have the following short time expansion.
In the following, $\operatorname{ev}_1$ denotes the evaluation map at time $1$,
that is, $\operatorname{ev}_1 (x)=x_1$ for
an ${\mathbf{R}}^n$-valued path $x$.
%
%co7.5 #&#
\begin{co}\label{coshort}
Let the coefficients $\sigma\dvtx  {\mathbf{R}}^n \to\operatorname{Mat} (n,d)$ and $\hat
\beta\dvtx  {\mathbf{R}}^n \to{\mathbf{R}}^n$
be $C_b^{\infty}$
and consider the RDE (\ref{eqrdefrac222}) above with $X= W^H$, where
$H \in(1/4,1/3) \cup(1/3,1/2)$.
Let $f$ and $g$ be real-valued $C_b^{\infty}$-functions on ${\mathbf{R}}^n$
such that
$F := f \circ\operatorname{ev}_1$ and $G := g \circ\operatorname{ev}_1$ satisfy
Assumptions \textup{(H1)--(H4)}.
Then,
we have the following asymptotic expansion as $t \searrow0$:
there are real constants $c$ and $\hat\alpha_{\kappa_0} (= \hat
\alpha
_0), \hat\alpha_{\kappa_1}, \hat\alpha_{\kappa_2},\ldots$ such that
\begin{eqnarray*}
&&
{\mathbb E} \bigl[
g( V^{ 1}_{0,t }) \exp\bigl( - f ( V^{ 1}_{0,t } ) / t^{2H} \bigr)
\bigr]
\\
&&\qquad=
\exp\bigl( - F_{\Lambda} (\gamma) /t^{2H} \bigr) \exp(- c/ t^H)
\\
&&\quad\qquad{}\times
\bigl(
\hat\alpha_{\kappa_0} +\hat\alpha_{\kappa_1} t^{\kappa_1 H} +
\cdots+
\hat\alpha_{\kappa_m} t^{\kappa_m H} +O( t^{\kappa_{m+1}H })
\bigr)
\end{eqnarray*}
for any $m \ge0$.
\end{co}

%re7.6 #&#
\begin{re}
Very roughly speaking,
in~\cite{bc,nnrt}, they studied the sort time asymptotics of the
following quantity under mild assumptions:
\[
{\mathbb E} [
g( V^{ 1}_{0,t } ) ].
\]
If $f$ is identically zero in Corollary~\ref{coshort},
then it is the same short time problems studied in~\cite{bc,nnrt}, at
least formally.
(It does not seem to the author that either~\cite{bc,nnrt} or the
Corollary~\ref{coshort} implies the other.)
\end{re}

%%%%%%%%%%%%%%%%%%%%%%%%%%%%%%%%%%%%%%%%%%%%%
%%%%%%%%%%%%%%%%%%%%%%%%%%%%%%%%%%%%%%%%%%%%%

%%%%%%%%%%%%%%%%%%%%%%%%%%%%%%%%%%%%%%%%%%%

% imsref loaded by akundreckaite, 2012-05-21 15:21:47

%suskaldyti doi

\printaddresses


\begin{thebibliography}{41}
% BibTex style file: ims.bst, 2011-05-30
% Default style options (sort=0,type=number).
% Used options (sort=1,type=number).

%b1 #&#
\bibitem{ad}
\begin{bbook}[mr]
\bauthor{\bsnm{Adams},~\bfnm{Robert~A.}\binits{R.~A.}}
(\byear{1975}).
\btitle{Sobolev Spaces}.
\bpublisher{Academic Press}, \baddress{New York--London}.
\bid{mr={0450957}}
\bptok{imsref}%
\end{bbook}
\endbibitem

%b2 #&#
\bibitem{aida}
\begin{barticle}[mr]
\bauthor{\bsnm{Aida},~\bfnm{Shigeki}\binits{S.}}
(\byear{2007}).
\btitle{Semi-classical limit of the bottom of spectrum of a {S}chr\"odinger
  operator on a path space over a compact {R}iemannian manifold}.
\bjournal{J. Funct. Anal.}
\bvolume{251}
\bpages{59--121}.
\bid{doi={10.1016/j.jfa.2007.06.009}, issn={0022-1236}, mr={2353701}}
\bptok{imsref}%
\end{barticle}
\endbibitem

%b3 #&#
\bibitem{ars}
\begin{barticle}[mr]
\bauthor{\bsnm{Albeverio},~\bfnm{Sergio}\binits{S.}},
  \bauthor{\bsnm{R{\"o}ckle},~\bfnm{Haio}\binits{H.}} \AND
  \bauthor{\bsnm{Steblovskaya},~\bfnm{Victoria}\binits{V.}}
(\byear{2000}).
\btitle{Asymptotic expansions for {O}rnstein--{U}hlenbeck semigroups perturbed
  by potentials over {B}anach spaces}.
\bjournal{Stochastics Stochastics Rep.}
\bvolume{69}
\bpages{195--238}.
\bid{issn={1045-1129}, mr={1760977}}
\bptok{imsref}%
\end{barticle}
\endbibitem

%b4 #&#
\bibitem{az}
\begin{bincollection}[mr]
\bauthor{\bsnm{Azencott},~\bfnm{Robert}\binits{R.}}
(\byear{1982}).
\btitle{Formule de {T}aylor stochastique et d\'eveloppement asymptotique
  d'int\'egrales de {F}eynman}.
In \bbooktitle{Seminar on {P}robability, {XVI}, {S}upplement}.
\bseries{Lecture Notes in Math.}
\bvolume{921}
\bpages{237--285}.
\bpublisher{Springer}, \baddress{Berlin}.
\bid{mr={0658728}}
\bptok{imsref}%
\end{bincollection}
\endbibitem

%b5 #&#
\bibitem{bc}
\begin{barticle}[mr]
\bauthor{\bsnm{Baudoin},~\bfnm{Fabrice}\binits{F.}} \AND
  \bauthor{\bsnm{Coutin},~\bfnm{Laure}\binits{L.}}
(\byear{2007}).
\btitle{Operators associated with a stochastic differential equation driven by
  fractional {B}rownian motions}.
\bjournal{Stochastic Process. Appl.}
\bvolume{117}
\bpages{550--574}.
\bid{doi={10.1016/j.spa.2006.09.004}, issn={0304-4149}, mr={2320949}}
\bptok{imsref}%
\end{barticle}
\endbibitem

%b6 #&#
\bibitem{bc2}
\begin{barticle}[mr]
\bauthor{\bsnm{Baudoin},~\bfnm{Fabrice}\binits{F.}} \AND
  \bauthor{\bsnm{Coutin},~\bfnm{Laure}\binits{L.}}
(\byear{2008}).
\btitle{Self-similarity and fractional {B}rownian motions on {L}ie groups}.
\bjournal{Electron. J. Probab.}
\bvolume{13}
\bpages{1120--1139}.
\bid{doi={10.1214/EJP.v13-530}, issn={1083-6489}, mr={2424989}}
\bptok{imsref}%
\end{barticle}
\endbibitem

%b7 #&#
\bibitem{ba}
\begin{barticle}[mr]
\bauthor{\bsnm{Ben~Arous},~\bfnm{G{\'e}rard}\binits{G.}}
(\byear{1988}).
\btitle{Methods de {L}aplace et de la phase stationnaire sur l'espace de
  {W}iener}.
\bjournal{Stochastics}
\bvolume{25}
\bpages{125--153}.
\bid{doi={10.1080/17442508808833536}, issn={0090-9491}, mr={0999365}}
\bptok{imsref}%
\end{barticle}
\endbibitem

%b8 #&#
\bibitem{bhoz}
\begin{bbook}[mr]
\bauthor{\bsnm{Biagini},~\bfnm{Francesca}\binits{F.}},
  \bauthor{\bsnm{Hu},~\bfnm{Yaozhong}\binits{Y.}},
  \bauthor{\bsnm{{\O}ksendal},~\bfnm{Bernt}\binits{B.}} \AND
  \bauthor{\bsnm{Zhang},~\bfnm{Tusheng}\binits{T.}}
(\byear{2008}).
\btitle{Stochastic Calculus for Fractional {B}rownian Motion and Applications}.
\bpublisher{Springer}, \baddress{London}.
\bid{doi={10.1007/978-1-84628-797-8}, mr={2387368}}
\bptok{imsref}%
\end{bbook}
\endbibitem

%b9 #&#
\bibitem{cou}
\begin{bincollection}[mr]
\bauthor{\bsnm{Coutin},~\bfnm{Laure}\binits{L.}}
(\byear{2007}).
\btitle{An introduction to (stochastic) calculus with respect to fractional
  {B}rownian motion}.
In \bbooktitle{S\'eminaire de {P}robabilit\'es {XL}}.
\bseries{Lecture Notes in Math.}
\bvolume{1899}
\bpages{3--65}.
\bpublisher{Springer}, \baddress{Berlin}.
\bid{doi={10.1007/978-3-540-71189-6_1}, mr={2408998}}
\bptok{imsref}%
\end{bincollection}
\endbibitem

%b10 #&#
\bibitem{cq}
\begin{barticle}[mr]
\bauthor{\bsnm{Coutin},~\bfnm{Laure}\binits{L.}} \AND
  \bauthor{\bsnm{Qian},~\bfnm{Zhongmin}\binits{Z.}}
(\byear{2002}).
\btitle{Stochastic analysis, rough path analysis and fractional {B}rownian
  motions}.
\bjournal{Probab. Theory Related Fields}
\bvolume{122}
\bpages{108--140}.
\bid{doi={10.1007/s004400100158}, issn={0178-8051}, mr={1883719}}
\bptok{imsref}%
\end{barticle}
\endbibitem

%%b11 #&#
%  \bauthor{\bsnm{{\"U}st{\"u}nel},~\bfnm{A.~S.}\binits{A.~S.}}
%(\byear{1999}).

%b12 #&#
\bibitem{eld}
\begin{bmisc}[auto:STB|2012/04/30|08:06:40]
\bauthor{\bsnm{Eldredge},~\bfnm{N.}\binits{N.}}
(\byear{2005}).
\bhowpublished{Computing $p$-variation. Unpublished manuscript. Univ.
  California, San Diego}.
\bptok{imsref}%
\end{bmisc}
\endbibitem

%b13 #&#
\bibitem{fo}
\begin{barticle}[mr]
\bauthor{\bsnm{Friz},~\bfnm{Peter}\binits{P.}} \AND
  \bauthor{\bsnm{Oberhauser},~\bfnm{Harald}\binits{H.}}
(\byear{2010}).
\btitle{A generalized {F}ernique theorem and applications}.
\bjournal{Proc. Amer. Math. Soc.}
\bvolume{138}
\bpages{3679--3688}.
\bid{doi={10.1090/S0002-9939-2010-10528-2}, issn={0002-9939}, mr={2661566}}
\bptok{imsref}%
\end{barticle}
\endbibitem

%b14 #&#
\bibitem{fv}
\begin{barticle}[mr]
\bauthor{\bsnm{Friz},~\bfnm{Peter}\binits{P.}} \AND
  \bauthor{\bsnm{Victoir},~\bfnm{Nicolas}\binits{N.}}
(\byear{2006}).
\btitle{A variation embedding theorem and applications}.
\bjournal{J. Funct. Anal.}
\bvolume{239}
\bpages{631--637}.
\bid{doi={10.1016/j.jfa.2005.12.021}, issn={0022-1236}, mr={2261341}}
\bptok{imsref}%
\end{barticle}
\endbibitem

%b15 #&#
\bibitem{fv0}
\begin{barticle}[auto:STB|2012/04/30|08:06:40]
\bauthor{\bsnm{Friz},~\bfnm{P.}\binits{P.}} \AND
  \bauthor{\bsnm{Victoir},~\bfnm{N.}\binits{N.}}
(\byear{2007}).
\btitle{Large deviation principle for enhanced Gaussian processes}.
\bjournal{Ann. l'Inst. Henri Poincar\'{e} Probab. Stat.}
\bvolume{43}
\bpages{775--785}.
\bptok{imsref}%
\end{barticle}
\endbibitem

%b16 #&#
\bibitem{fvgauss}
\begin{barticle}[mr]
\bauthor{\bsnm{Friz},~\bfnm{Peter}\binits{P.}} \AND
  \bauthor{\bsnm{Victoir},~\bfnm{Nicolas}\binits{N.}}
(\byear{2010}).
\btitle{Differential equations driven by {G}aussian signals}.
\bjournal{Ann. Inst. Henri Poincar\'e Probab. Stat.}
\bvolume{46}
\bpages{369--413}.
\bid{doi={10.1214/09-AIHP202}, issn={0246-0203}, mr={2667703}}
\bptok{imsref}%
\end{barticle}
\endbibitem

%b17 #&#
\bibitem{fvbook}
\begin{bbook}[mr]
\bauthor{\bsnm{Friz},~\bfnm{Peter~K.}\binits{P.~K.}} \AND
  \bauthor{\bsnm{Victoir},~\bfnm{Nicolas~B.}\binits{N.~B.}}
(\byear{2010}).
\btitle{Multidimensional Stochastic Processes as Rough Paths: Theory and Applications}.
\bseries{Cambridge Studies in Advanced Mathematics}
\bvolume{120}.
\bpublisher{Cambridge Univ. Press}, \baddress{Cambridge}.
\bid{mr={2604669}}
\bptok{imsref}%
\end{bbook}
\endbibitem

%b18 #&#
\bibitem{glj}
\begin{bmisc}[auto:STB|2012/04/30|08:06:40]
\bauthor{\bsnm{Gubinelli},~\bfnm{M.}\binits{M.}} \AND
  \bauthor{\bsnm{Lejay},~\bfnm{A.}\binits{A.}}
  (\byear{2009}).
\bhowpublished{Global existence for rough differential equations under linear
  growth condition. Preprint}.
\bptok{imsref}%
\end{bmisc}
\endbibitem

%b19 #&#
\bibitem{ina}
\begin{barticle}[mr]
\bauthor{\bsnm{Inahama},~\bfnm{Yuzuru}\binits{Y.}}
(\byear{2006}).
\btitle{Laplace's method for the laws of heat processes on loop spaces}.
\bjournal{J. Funct. Anal.}
\bvolume{232}
\bpages{148--194}.
\bid{doi={10.1016/j.jfa.2005.06.006}, issn={0022-1236}, mr={2200170}}
\bptok{imsref}%
\end{barticle}
\endbibitem

%b20 #&#
\bibitem{ina2}
\begin{barticle}[mr]
\bauthor{\bsnm{Inahama},~\bfnm{Yuzuru}\binits{Y.}}
(\byear{2010}).
\btitle{A stochastic {T}aylor-like expansion in the rough path theory}.
\bjournal{J. Theoret. Probab.}
\bvolume{23}
\bpages{671--714}.
\bid{doi={10.1007/s10959-010-0287-6}, issn={0894-9840}, mr={2679952}}
\bptok{imsref}%
\end{barticle}
\endbibitem

%b21 #&#
\bibitem{ik}
\begin{barticle}[mr]
\bauthor{\bsnm{Inahama},~\bfnm{Yuzuru}\binits{Y.}} \AND
  \bauthor{\bsnm{Kawabi},~\bfnm{Hiroshi}\binits{H.}}
(\byear{2007}).
\btitle{Asymptotic expansions for the {L}aplace approximations for {I}t\^o
  functionals of {B}rownian rough paths}.
\bjournal{J. Funct. Anal.}
\bvolume{243}
\bpages{270--322}.
\bid{doi={10.1016/j.jfa.2006.09.016}, issn={0022-1236}, mr={2291439}}
\bptok{imsref}%
\end{barticle}
\endbibitem

%b22 #&#
\bibitem{jain}
\begin{barticle}[mr]
\bauthor{\bsnm{Jain},~\bfnm{Naresh~C.}\binits{N.~C.}} \AND
  \bauthor{\bsnm{Monrad},~\bfnm{Ditlev}\binits{D.}}
(\byear{1983}).
\btitle{Gaussian measures in {$B\sb{p}$}}.
\bjournal{Ann. Probab.}
\bvolume{11}
\bpages{46--57}.
\bid{issn={0091-1798}, mr={0682800}}
\bptok{imsref}%
\end{barticle}
\endbibitem

%b23 #&#
\bibitem{jan}
\begin{bbook}[mr]
\bauthor{\bsnm{Janson},~\bfnm{Svante}\binits{S.}}
(\byear{1997}).
\btitle{Gaussian {H}ilbert Spaces}.
\bseries{Cambridge Tracts in Mathematics}
\bvolume{129}.
\bpublisher{Cambridge Univ. Press}, \baddress{Cambridge}.
\bid{doi={10.1017/CBO9780511526169}, mr={1474726}}
\bptok{imsref}%
\end{bbook}
\endbibitem

%b24 #&#
\bibitem{kuo}
\begin{bbook}[mr]
\bauthor{\bsnm{Kuo},~\bfnm{Hui~Hsiung}\binits{H.~H.}}
(\byear{1975}).
\btitle{Gaussian Measures in {B}anach Spaces}.
\bseries{Lecture Notes in Math.}
\bvolume{463}.
\bpublisher{Springer}, \baddress{Berlin}.
\bid{mr={0461643}}
\bptok{imsref}%
\end{bbook}
\endbibitem

%b25 #&#
\bibitem{ko}
\begin{barticle}[mr]
\bauthor{\bsnm{Kusuoka},~\bfnm{Shigeo}\binits{S.}} \AND
  \bauthor{\bsnm{Osajima},~\bfnm{Yasufumi}\binits{Y.}}
(\byear{2008}).
\btitle{A remark on the asymptotic expansion of density function of {W}iener
  functionals}.
\bjournal{J. Funct. Anal.}
\bvolume{255}
\bpages{2545--2562}.
\bid{doi={10.1016/j.jfa.2008.03.019}, issn={0022-1236}, mr={2473267}}
\bptok{imsref}%
\end{barticle}
\endbibitem

%b26 #&#
\bibitem{ks1}
\begin{barticle}[mr]
\bauthor{\bsnm{Kusuoka},~\bfnm{Shigeo}\binits{S.}} \AND
  \bauthor{\bsnm{Stroock},~\bfnm{Daniel~W.}\binits{D.~W.}}
(\byear{1991}).
\btitle{Precise asymptotics of certain {W}iener functionals}.
\bjournal{J. Funct. Anal.}
\bvolume{99}
\bpages{1--74}.
\bid{doi={10.1016/0022-1236(91)90051-6}, issn={0022-1236}, mr={1120913}}
\bptok{imsref}%
\end{barticle}
\endbibitem

%b27 #&#
\bibitem{ks2}
\begin{barticle}[mr]
\bauthor{\bsnm{Kusuoka},~\bfnm{Shigeo}\binits{S.}} \AND
  \bauthor{\bsnm{Stroock},~\bfnm{Daniel~W.}\binits{D.~W.}}
(\byear{1994}).
\btitle{Asymptotics of certain {W}iener functionals with degenerate extrema}.
\bjournal{Comm. Pure Appl. Math.}
\bvolume{47}
\bpages{477--501}.
\bid{doi={10.1002/cpa.3160470404}, issn={0010-3640}, mr={1272385}}
\bptok{imsref}%
\end{barticle}
\endbibitem

%%b28 #&#
%  \bauthor{\bsnm{Qian},~\bfnm{Z.}\binits{Z.}} \AND
%  \bauthor{\bsnm{Zhang},~\bfnm{T.}\binits{T.}}
%(\byear{2002}).
%  paths}.

%b29 #&#
\bibitem{lej}
\begin{bincollection}[mr]
\bauthor{\bsnm{Lejay},~\bfnm{Antoine}\binits{A.}}
(\byear{2003}).
\btitle{An introduction to rough paths}.
In \bbooktitle{S\'eminaire de {P}robabilit\'es {XXXVII}}.
\bseries{Lecture Notes in Math.}
\bvolume{1832}
\bpages{1--59}.
\bpublisher{Springer}, \baddress{Berlin}.
\bid{doi={10.1007/978-3-540-40004-2_1}, mr={2053040}}
\bptok{imsref}%
\end{bincollection}
\endbibitem

%b30 #&#
\bibitem{ll}
\begin{barticle}[mr]
\bauthor{\bsnm{Li},~\bfnm{Xiang-Dong}\binits{X.-D.}} \AND
  \bauthor{\bsnm{Lyons},~\bfnm{Terry~J.}\binits{T.~J.}}
(\byear{2006}).
\btitle{Smoothness of {I}t\^o maps and diffusion processes on path spaces.
  {I}}.
\bjournal{Ann. Sci. \'Ec. Norm. Sup\'{e}r. (4)}
\bvolume{39}
\bpages{649--677}.
\bid{doi={10.1016/j.ansens.2006.07.001}, issn={0012-9593}, mr={2290140}}
\bptok{imsref}%
\end{barticle}
\endbibitem

%b31 #&#
\bibitem{lq}
\begin{bbook}[mr]
\bauthor{\bsnm{Lyons},~\bfnm{Terry}\binits{T.}} \AND
  \bauthor{\bsnm{Qian},~\bfnm{Zhongmin}\binits{Z.}}
(\byear{2002}).
\btitle{System Control and Rough Paths}.
\bpublisher{Oxford Univ. Press}, \baddress{Oxford}.
\bid{doi={10.1093/acprof:oso/9780198506485.001.0001}, mr={2036784}}
\bptok{imsref}%
\end{bbook}
\endbibitem

%b32 #&#
\bibitem{ly}
\begin{barticle}[mr]
\bauthor{\bsnm{Lyons},~\bfnm{Terry~J.}\binits{T.~J.}}
(\byear{1998}).
\btitle{Differential equations driven by rough signals}.
\bjournal{Rev. Mat. Iberoam.}
\bvolume{14}
\bpages{215--310}.
\bid{doi={10.4171/RMI/240}, issn={0213-2230}, mr={1654527}}
\bptok{imsref}%
\end{barticle}
\endbibitem

%b33 #&#
\bibitem{lcl}
\begin{bbook}[mr]
\bauthor{\bsnm{Lyons},~\bfnm{Terry~J.}\binits{T.~J.}},
  \bauthor{\bsnm{Caruana},~\bfnm{Michael}\binits{M.}} \AND
  \bauthor{\bsnm{L{\'e}vy},~\bfnm{Thierry}\binits{T.}}
(\byear{2007}).
\btitle{Differential Equations Driven by Rough Paths}.
\bseries{Lecture Notes in Math.}
\bvolume{1908}.
\bpublisher{Springer}, \baddress{Berlin}.
%  Saint-Flour, July 6--24, 2004, With an introduction concerning the Summer
%  School by Jean Picard}.
\bid{mr={2314753}}
\bptok{imsref}%
\end{bbook}
\endbibitem

%b34 #&#
\bibitem{mss}
\begin{barticle}[mr]
\bauthor{\bsnm{Millet},~\bfnm{Annie}\binits{A.}} \AND
  \bauthor{\bsnm{Sanz-Sol{\'e}},~\bfnm{Marta}\binits{M.}}
(\byear{2006}).
\btitle{Large deviations for rough paths of the fractional {B}rownian motion}.
\bjournal{Ann. Inst. Henri Poincar\'e Probab. Stat.}
\bvolume{42}
\bpages{245--271}.
\bid{doi={10.1016/j.anihpb.2005.04.003}, issn={0246-0203}, mr={2199801}}
\bptok{imsref}%
\end{barticle}
\endbibitem

%b35 #&#
\bibitem{mis}
\begin{bbook}[mr]
\bauthor{\bsnm{Mishura},~\bfnm{Yuliya~S.}\binits{Y.~S.}}
(\byear{2008}).
\btitle{Stochastic Calculus for Fractional {B}rownian Motion and Related
  Processes}.
\bseries{Lecture Notes in Math.}
\bvolume{1929}.
\bpublisher{Springer}, \baddress{Berlin}.
\bid{doi={10.1007/978-3-540-75873-0}, mr={2378138}}
\bptok{imsref}%
\end{bbook}
\endbibitem

%b36 #&#
\bibitem{nnrt}
\begin{barticle}[mr]
\bauthor{\bsnm{Neuenkirch},~\bfnm{A.}\binits{A.}},
  \bauthor{\bsnm{Nourdin},~\bfnm{I.}\binits{I.}},
  \bauthor{\bsnm{R{\"o}{\ss}ler},~\bfnm{A.}\binits{A.}} \AND
  \bauthor{\bsnm{Tindel},~\bfnm{S.}\binits{S.}}
(\byear{2009}).
\btitle{Trees and asymptotic expansions for fractional stochastic differential
  equations}.
\bjournal{Ann. Inst. Henri Poincar\'e Probab. Stat.}
\bvolume{45}
\bpages{157--174}.
\bid{doi={10.1214/07-AIHP159}, issn={0246-0203}, mr={2500233}}
\bptok{imsref}%
\end{barticle}
\endbibitem

%%b37 #&#
%(\byear{1995}).
%%\bseries{Probability and Its Applications (New York)}.

%b38 #&#
\bibitem{rt}
\begin{barticle}[mr]
\bauthor{\bsnm{Rovira},~\bfnm{Carles}\binits{C.}} \AND
  \bauthor{\bsnm{Tindel},~\bfnm{Samy}\binits{S.}}
(\byear{2000}).
\btitle{Sharp {L}aplace asymptotics for a parabolic {SPDE}}.
\bjournal{Stochastics Stochastics Rep.}
\bvolume{69}
\bpages{11--30}.
\bid{issn={1045-1129}, mr={1751716}}
\bptok{imsref}%
\end{barticle}
\endbibitem

%b39 #&#
\bibitem{tw}
\begin{bincollection}[mr]
\bauthor{\bsnm{Takanobu},~\bfnm{S.}\binits{S.}} \AND
  \bauthor{\bsnm{Watanabe},~\bfnm{S.}\binits{S.}}
(\byear{1993}).
\btitle{Asymptotic expansion formulas of the {S}childer type for a class of
  conditional {W}iener functional integrations}.
In \bbooktitle{Asymptotic Problems in Probability Theory: {W}iener Functionals
  and Asymptotics ({S}anda/{K}yoto, 1990)}.
\bseries{Pitman Res. Notes Math. Ser.}
\bvolume{284}
\bpages{194--241}.
\bpublisher{Longman Sci. Tech.}, \baddress{Harlow}.
\bid{mr={1354169}}
\bptok{imsref}%
\end{bincollection}
\endbibitem

%b40 #&#
\bibitem{tan}
\begin{barticle}[mr]
\bauthor{\bsnm{Taniguchi},~\bfnm{Setsuo}\binits{S.}}
(\byear{2008}).
\btitle{Quadratic {W}iener functionals of square norms on measure spaces}.
\bjournal{Commun. Stoch. Anal.}
\bvolume{2}
\bpages{11--26}.
\bid{issn={0973-9599}, mr={2446908}}
\bptok{imsref}%
\end{barticle}
\endbibitem

\end{thebibliography}
\end{document}